# COHESIVE FRACTURE WITH IRREVERSIBILITY:
# QUASISTATIC EVOLUTION FOR A MODEL SUBJECT TO FATIGUE

VITO CRISMALE, GIULIANO LAZZARONI, AND GIANLUCA ORLANDO

ABSTRACT. In this paper we prove the existence of quasistatic evolutions for a cohesive fracture on a prescribed crack surface, in small-strain antiplane elasticity. The main feature of the model is that the density of the energy dissipated in the fracture process depends on the total variation of the amplitude of the jump. Thus, any change in the crack opening entails a loss of energy, until the crack is complete. In particular this implies a fatigue phenomenon, i.e., a complete fracture may be produced by oscillation of small jumps.

The first step of the existence proof is the construction of approximate evolutions obtained by solving discrete-time incremental minimum problems. The main difficulty in the passage to the continuous-time limit is that we lack of controls on the variations of the jump of the approximate evolutions. Therefore we resort to a weak formulation where the variation of the jump is replaced by a Young measure. Eventually, after proving the existence in this weak formulation, we improve the result by showing that the Young measure is concentrated on a function and coincides with the variation of the jump of the displacement.

## CONTENTS

## 1. INTRODUCTION

In this paper we investigate on the quasistatic evolution of cohesive cracks in a material subject to a fatigue phenomenon. Compared to brittle fracture, cohesive models provide a more accurate description of the process of crack growth. Indeed, in Griffith's theory of brittle fracture [24], the energy spent to produce a crack only depends on the geometry of the crack itself, the simplest case being a surface energy proportional to the measure of the crack set. In contrast, cohesive energies, introduced in [7], also depend on the crack opening, i.e., on the difference between the traces of the displacement on the two sides of the crack. In fact, fracture should be regarded as a gradual process, where the material is considered completely cracked at a point only when the amplitude of the jump of the displacement is sufficiently large.

Thus, when the crack opening gradually increases in time, some energy is dissipated, until the opening overcomes a certain threshold. On the other hand, mechanical systems may present different responses when the crack opening happens to decrease: in our model, some energy is dissipated also in this phase (until a maximal dissipation is reached), because of the contact between the two sides of the crack. In this respect, the behaviour of our system differs from those considered in the mathematical literature on quasistatic cohesive fracture. For instance, when the crack opening decreases one may assume that no energy is dissipated [21] or that some dissipated energy is recovered [11, 4]. (See also Remark 2.5.) The peculiar response modelled here leads to a fatigue phenomenon affecting crack growth and gives further mathematical difficulties as we outline below.

We now describe the formulation of the problem in the setting of small-strain antiplane elasticity, referring to Section 2 for all details on the mathematical assumptions. More precisely, the reference configuration of the body is supposed to be an infinite cylinder $\Omega \times \mathbb{R}$, with $\Omega \subset \mathbb{R}^n$ ($n = 2$ being the physically relevant case), and the deformation $v \colon \Omega \times \mathbb{R} \to \mathbb{R}^{n+1}$ takes the form $v(x_1, \ldots, x_{n+1}) = (x_1, \ldots, x_n, x_{n+1} + u(x_1, \ldots, x_n))$, where $u \colon \Omega \to \mathbb{R}$ is the vertical displacement. According to the small-strain assumption, the linearized elastic energy is

$$\frac{1}{2} \int_\Omega |\nabla u|^2 \, \mathrm{d}x \,.$$

The body may present cracks of the form $\hat{\Gamma} \times \mathbb{R}$, where $\hat{\Gamma}$ is contained in a prescribed $(n-1)$-dimensional manifold $\Gamma \subset \mathbb{R}^n$. Removing the restriction of a prescribed fracture set is by now out of reach in the mathematical modelling of cohesive crack growth, in contrast to the brittle case, where several existence results were obtained under more general hypotheses [20, 12, 22, 17, 18, 26].

As mentioned above, the energy dissipated during the fracture process here depends on the evolution of the amplitude of the jump, denoted by $[u(t)] \colon \Gamma \to \mathbb{R}$, where $t \in [0, T]$ is the time variable. To describe the response of the system to loading, we start by considering the situation where $[u(0)] = 0$ on $\Gamma$ and $t \mapsto [u(t)]$ is nondecreasing on $\Gamma$ in a time interval $[0, t_1]$. In this case, the energy dissipated in $[0, t_1]$ is

$$\int_\Gamma g\big(\big|[u(t_1)]\big|\big) \, \mathrm{d}\mathcal{H}^{n-1},$$

where $g \colon [0, +\infty) \to [0, +\infty)$ is a concave (thus nondecreasing) function satisfying: $g(0) = 0$; $g'(0)$ exists, finite; and $g(\xi) \to \kappa \in (0, \infty)$ as $\xi \to \infty$. (See Section 2 for more general assumptions on $g$.) If, afterwards, $t \to [u(t)]$ is nonincreasing in the interval $[t_1, t_2]$, there is still some dissipated energy in $[t_1, t_2]$, which amounts to

$$\int_\Gamma g\big(\big|[u(t_1)]\big| + \big|[u(t_2)] - [u(t_1)]\big|\big) \, \mathrm{d}\mathcal{H}^{n-1} - \int_\Gamma g\big(\big|[u(t_1)]\big|\big) \, \mathrm{d}\mathcal{H}^{n-1}.$$

As a consequence, a complete fracture (corresponding to $g = \kappa$) may occur not only after a large crack opening, but even after oscillations of small jumps (e.g. by a cyclic loading).

In fact, on the contact area between the two cracked parts of the material, the repeated relative surface motion can induce damage by a fatigue process. In applications, this wear phenomenon is known as fretting [14] and occurs as a result of relative sliding motion of the order from nanometres to millimetres.

Fatigue effects related to cohesive fracture have been already observed e.g. in [2, 3], where a cohesive crack appears in an elastoplastic material subject to damage. In that model, damage occurs more easily in regions where the material has suffered cyclic plastic deformations. Moreover, a cohesive law relates damage with the amplitude of the jump of the displacement (i.e., the plastic slip). Thus, due to the irreversibility of damage, oscillations of the crack opening result in energy dissipation until damage is complete. The existence of quasistatic evolution for coupled elastoplastic-damage models has been proved in [15, 16]. The limit of such a model when damage is forced to concentrate on hypersurfaces has been studied in [19] in a static setting and gives rise to cohesive surface energies. We also refer to [13] and [6] for the derivation of cohesive-type energies by means of $\Gamma$-convergence.

These motivations lead us to formulate a discrete-time problem as already proposed in [1]. Given an initial condition $u(0) = u_0$ and a time-dependent Dirichlet datum $w(t)$ on $\partial_D \Omega \subset \partial \Omega$ (cf. Section 2 for the complete list of assumptions), for every $k \in \mathbb{N}$ we fix a subdivision $0 = t_k^0 < t_k^1 <$

$\cdots < t_k^{k-1} < t_k^k = T$ and we define recursively $u_k^i$ and $V_k^i$ by

$$u_k^i \in \operatorname*{argmin}_u \left\{ \frac{1}{2} \int_\Omega |\nabla u|^2 \, \mathrm{d}x + \int_\Gamma g\Big(V_k^{i-1} + \big|[u] - [u_k^{i-1}]\big|\Big) \, \mathrm{d}\mathcal{H}^{n-1} \ : \ u = w(t_k^i) \text{ on } \partial_D \Omega \right\}, \quad (1.1)$$

$$V_k^i := V_k^{i-1} + \big|[u_k^i] - [u_k^{i-1}]\big|,$$

where $u_k^0 := u_0$ and $V_k^0 = |[u_0]|$. The function $V_k^i$ describes the cumulated jump of the approximate evolutions at each point of $\Gamma$. The authors wish to thank Jean–Jacques Marigo for pointing out that the model studied in this paper is connected with the discrete model proposed in [1].

We shall pass to the limit as $k \to \infty$ and prove that the resulting continuous-time evolution satisfies the usual properties of quasistatic processes: stability and energy balance. This strategy has been common when proving existence of quasistatic evolutions in fracture mechanics since the seminal paper [23] (see also [8]), and, more in general, when looking for energetic solutions to rate-independent systems [27]. Loosely speaking, following this approach one selects equilibria of the system among the global minimisers of the sum of the mechanical energy and of the dissipated energy. Such restriction is advantageous from the mathematical point of view, but may lead to unphysical phenomena when the solution presents time discontinuities. For this reason one may resort to notions of solutions based on local minimality, e.g. by means of a vanishing viscosity approach, cf. [25, 26] in the brittle case and [10, 4] for a cohesive model.

In order to study the continuous-time limit, we define $u_k(t)$ and $V_k(t)$ as the piecewise constant interpolations of $u_k^i$ and $V_k^i$ in time, respectively. The main difficulty in the passage to limit as $k \to \infty$ is that we lack of controls on $V_k(t)$. In fact, by (1.1), we can only infer that $\int_\Gamma g(V_k(t)) \, \mathrm{d}\mathcal{H}^{n-1}$ is uniformly bounded, but this gives no information on the equi-integrability of $V_k(t)$, since $g$ is bounded. (This would not be the case if $g$ had e.g. linear growth as in a model for perfect plasticity constrained on $\Gamma$.) In the first instance, in order to pass to the limit as $k \to \infty$, the only chance is to employ compactness properties of the wider class of Young measures (as already done in [11]). Indeed, because of the monotonicity of $V_k(t)$, a Helly-type selection principle [11, Theorem 3.20] guarantees that $V_k(t)$ generates a Young measure $\nu(t) = (\nu^x(t))_{x\in\Gamma}$ for every $t$, up to a subsequence independent of $t$.

As for the displacements, from the uniform a priori bounds we obtain that there is a subsequence $u_{k_j}(t)$ weakly converging to a function $u(t)$. Yet the subsequence $k_j = k_j(t)$ may depend on $t$. This is a technical inconvenience, since we need to keep track of the relation between $V_k(t)$ and $[u_k(t)]$, namely, to pass to the limit in the irreversibility relation

$$V_k(t) \geq V_k(s) + \big|[u_k(t)] - [u_k(s)]\big| \quad \text{for any } s \leq t.$$

Indeed, notice that $u_k(t)$ and $u_k(s)$ may converge along different subsequences! This difficulty is solved by rewriting the previous inequality as a system of two inequalities

$$V_k(t) + [u_k(t)] \geq V_k(s) + [u_k(s)] \quad \text{for any } s \leq t, \qquad (1.2)$$

$$V_k(t) - [u_k(t)] \geq V_k(s) - [u_k(s)] \quad \text{for any } s \leq t. \qquad (1.3)$$

In fact, we can now pass to the limit in these relations by means of a Helly-type theorem, extracting a further subsequence (not relabelled) independent of $t$ and exploiting the monotonicity of $V_k(t) \pm [u_k(t)]$. Moreover, thanks to this trick it turns out that we can identify the limit jump $[u(t)]$ without extracting further subsequences. Ultimately, also the displacement $u(t)$ is the limit of the whole sequence $u_k(t)$, since $u(t)$ is the solution of a minimum problem among functions with prescribed jump $[u(t)]$. This property is relevant for the approximation of the solutions.

At this point of the analysis, we can pass to the limit in the global stability and in the energy balance, obtaining that $(u(t), \nu(t))$ fulfils a weak notion of quasistatic evolution. Specifically, the variation of jumps on $\Gamma$ is replaced by the Young measure $\nu(t)$. On the other hand, we can define the variation of the jumps of $u(t)$, denoted by $V_u(t)$. We refer to Section 2 for a rigorous definition; here we only notice that, if $u$ is an absolutely continuous function of time, then we have

$$V_u(t) = \int_0^t \big|[\dot{u}(s)]\big| \, \mathrm{d}s.$$

(See Lemma 6.2.) Actually, we improve the existence result by proving that $(u(t), V_u(t))$ satisfies the same properties of global stability and energy balance: this shows the existence of a quasistatic evolution in a stronger formulation that does not employ Young measures. Furthermore, we prove that (up to a truncation) the limit measure $\nu(t)$ is concentrated on the limit function $V_u(t)$, so also the limit of the discrete variations $V_k(t)$ is characterised.

The notion of quasistatic evolution and the main existence result are presented in Section 2, which contains also some results on a strong formulation that is satisfied by the weak solutions under suitable regularity assumptions. The final part of Section 2 contains a short presentation of the existence proof, which is given in more detail in the remaining part of the paper. After recalling some preliminary results on Young measures (Section 3), we introduce the discrete-time problems in Section 4 and we pass to the continuous-time limit in Section 5, obtaining the formulation based on Young measures. Finally, in Section 6 we prove the existence of quasistatic evolutions according to the notion based on functions.

## 2. Assumptions on the model and statement of the main result

**Notation.** If $\Xi$ is a metric space, we denote by $\mathscr{B}(\Xi)$ the $\sigma$-algebra of Borel sets on $\Xi$. The Lebesgue measure in $\mathbb{R}^n$ is denoted by $\mathcal{L}^n$, while $\mathcal{H}^{n-1}$ is the $(n-1)$-dimensional Hausdorff measure.

Given a Hilbert space $X$, we recall that $AC([0,T];X)$ is the space of all absolutely continuous functions defined in $[0,T]$ with values in $X$. For the main properties of these functions we refer, e.g., to [9, Appendix]. Given $\gamma \in AC([0,T];X)$, the time derivative of $\gamma$, which exists a.e. in $[0,T]$, is denoted by $\dot{\gamma}$. It is well-known that $\dot{\gamma}$ is a Bochner integrable function with values in $X$.

In the sequel, we will often consider time-dependent functions $t \mapsto v(t)$, where $v(t)$ is a function depending on a space variable $x$. We will write $v(t;x)$ to refer to the value of $v(t)$ in $x$.

**Reference configuration and boundary conditions.** Throughout the paper, $\Omega$ is a bounded, Lipschitz, open set in $\mathbb{R}^n$ representing the cross-section of a cylindrical body in the reference configuration (in the setting of antiplane shear). The cracks of the body will be contained in a prescribed crack surface $\Gamma$, where $\Gamma$ is a $(n-1)$-dimensional Lipschitz manifold in $\mathbb{R}^n$ with $0 < \mathcal{H}^{n-1}(\Gamma \cap \overline{\Omega}) < \infty$. Moreover, we assume that $\Omega \setminus \Gamma = \Omega^+ \cup \Omega^-$, where $\Omega^+$ and $\Omega^-$ are disjoint open connected sets with Lipschitz boundary. The normal $\nu(x) = \nu_\Gamma(x)$ to the surface $\Gamma$ is chosen in such a way that it coincides with the outer normal to $\partial\Omega^-$.

We consider evolutions driven by a time-dependent boundary condition assigned on the Dirichlet part of the boundary $\partial_D\Omega$. We assume that $\partial_D\Omega$ is a relatively open set of $\partial\Omega$ and that $\mathcal{H}^{n-1}(\partial_D\Omega \cap \partial\Omega^\pm) > 0$, in order to apply the Poincaré Inequality separately in $\Omega^+$ and $\Omega^-$. We denote by $\partial_N\Omega$ the remaining part of the boundary, i.e., $\partial_N\Omega := \partial\Omega \setminus \partial_D\Omega$.

For every $w \in H^1(\Omega)$, we define the set of *admissible displacements* corresponding to $w$ by

$$\mathscr{A}(w) := \{u \in H^1(\Omega \setminus \Gamma) \ : \ u = w \text{ on } \partial_D\Omega\}. \tag{2.1}$$

We assign a function $t \mapsto w(t)$ defined on $[0,T]$ with values in $H^1(\Omega)$ and we assume that

$$t \mapsto w(t) \quad \text{belongs to} \quad AC([0,T];H^1(\Omega)). \tag{2.2}$$

**Variation of jumps and initial data.** In order to give the notion of quasistatic evolution, we introduce a function $V_u(t)$ describing the variation of the jumps on $\Gamma$ of an evolution $t \mapsto u(t)$ in a time interval $[0,t]$.

Before defining $V_u(t)$, we recall the definition of the essential supremum of a family of measurable functions, that is the least upper bound in the sense of a.e. inequality. We give this definition in the case of functions defined on the measure space $(\Gamma; \mathcal{H}^{n-1})$. Indeed, this will be the relevant setting for our model.

**Definition 2.1.** Let $(v_i)_{i \in I}$ be a family of measurable functions from $\Gamma$ to $[-\infty, \infty]$. Let $\overline{v} \colon \Gamma \to [-\infty, \infty]$ be a measurable function such that

  (i) $\overline{v} \geq v_i$ $\mathcal{H}^{n-1}$-a.e. on $\Gamma$, for every $i \in I$;
  (ii) if $v \colon \Gamma \to [-\infty, \infty]$ is a measurable function such that $v \geq v_i$ $\mathcal{H}^{n-1}$-a.e. on $\Gamma$, for every $i \in I$, then $v \geq \overline{v}$ $\mathcal{H}^{n-1}$-a.e. on $\Gamma$.

We say that $\overline{v}$ an *essential supremum* of the family $(v_i)_{i \in I}$.

*Remark* 2.2. Given a family of measurable functions $(v_i)_{i \in I}$, there exists a unique (up to $\mathcal{H}^{n-1}$-a.e. equivalence) essential supremum $\overline{v}$ of the family $(v_i)_{i \in I}$. We denote it by $\operatorname*{ess\,sup}_{i \in I} v_i := \overline{v}$.

We now define the essential variation, namely the variation for a time-dependent family of measurable functions, in the sense of a.e. inequality. As done for the essential supremum, we give this definition in the case of functions defined on the measure space $(\Gamma; \mathcal{H}^{n-1})$.

**Definition 2.3.** Let us consider a function $t \mapsto \gamma(t)$, with $\gamma(t) \colon \Gamma \to \mathbb{R}$ measurable for every $t \in [0, T]$. For every $0 \leq t_1 \leq t_2 \leq T$, the *essential variation* of $\gamma$ in $[t_1, t_2]$ is the function $\operatorname{ess Var}(\gamma; t_1, t_2) \colon \Gamma \to [0, \infty]$ defined by

$$\operatorname{ess Var}(\gamma; t_1, t_2) := \operatorname{ess\,sup} \Big\{ \sum_{i=1}^{j} |\gamma(s_i) - \gamma(s_{i-1})| \ : \ j \in \mathbb{N}, \ t_1 = s_0 < s_1 < \cdots < s_{j-1} < s_j = t_2 \Big\}.$$

*Remark* 2.4. The essential variation satisfies the usual property that

$$\operatorname{ess Var}(\gamma; t_1, t_3) = \operatorname{ess Var}(\gamma; t_1, t_2) + \operatorname{ess Var}(\gamma; t_2, t_3) \quad \mathcal{H}^{n-1}\text{-a.e. on } \Gamma,$$

for any $0 \leq t_1 < t_2 < t_3 \leq t$.

Given a function $t \mapsto u(t)$ defined on $[0, T]$ with values in $H^1(\Omega \setminus \Gamma)$, we define the variation $V_u(t) \colon \Gamma \to [0, \infty]$ of its jumps on $\Gamma$ with initial condition $V_0$ by

$$V_u(t) := \operatorname{ess Var}([u]; 0, t) + V_0, \tag{2.3}$$

for every $t \in [0, T]$, where $V_0 \colon \Gamma \to [0, \infty]$ is an assigned measurable function.

**Initial data.** We fix an initial displacement

$$u_0 \in \mathscr{A}(w(0)) \tag{2.4}$$

and a function $V_0 \colon \Gamma \to [0, \infty]$ accounting for the variation of previous jumps until the initial time $t = 0$. Indeed we assume that

$$V_0(x) \geq \big|[u_0(x)]\big| \quad \text{for } \mathcal{H}^{n-1}\text{-a.e. } x \in \Gamma. \tag{2.5}$$

If $V_0 = \big|[u_0]\big|$, a monotone crack opening has occurred before the initial time $t = 0$. In general, the crack opening may have oscillated before the initial time in such a way that its variation in time equals $V_0$. The set $\Gamma_N(0) := \{V_0 \geq \theta(x)\}$ represents the part of $\Gamma$ which is already completely broken at the beginning of the process.

**The surface energy density.** We assume that the surface energy density $g$ depends on the point on $\Gamma$ and on the history of the jump. More precisely, $g \colon \Gamma \times [0, \infty) \to [0, \infty)$ satisfies the following assumptions:

$(g1)$ $g$ is a Carathéodory integrand, i.e., $g(x, \cdot)$ is continuous for $\mathcal{H}^{n-1}$-a.e. $x \in \Gamma$ and $g(\cdot, \xi)$ is $\mathcal{H}^{n-1}$-measurable for every $\xi \in [0, \infty)$;

$(g2)$ $g(x, 0) = 0$ and $g(x, \cdot)$ is concave for $\mathcal{H}^{n-1}$-a.e. $x \in \Gamma$;

$(g3)$ $\lim_{\xi \to \infty} g(x, \xi) = \kappa(x) \in [\kappa_1, \kappa_2]$ for $\mathcal{H}^{n-1}$-a.e. $x \in \Gamma$, where $\kappa_1, \kappa_2 \in (0, \infty)$;

$(g4)$ the limit

$$\lim_{\xi \to 0^+} \frac{g(x, \xi)}{\xi} =: g'(x, 0)$$

exists for $\mathcal{H}^{n-1}$-a.e. $x \in \Gamma$ and $g'(\cdot, 0) \in L^\infty(\Gamma)$.

In particular, for $\mathcal{H}^{n-1}$-a.e. $x \in \Gamma$ it turns out that $g(x, \cdot)$ is nondecreasing and can be extended to a function in $\mathcal{C}_b([0, \infty])$ by setting $g(x, \infty) := \kappa(x)$.

It will be convenient to introduce a measurable function $\theta \colon \Gamma \to [0, \infty]$ that represents the threshold after which the function $g(x, \cdot)$ becomes constant, i.e.,

$$\theta(x) := \inf\{\xi > 0 \ : \ g(x, \xi) = \kappa(x)\} \in (0, \infty]. \tag{2.6}$$

The function $g(x, \cdot)$ is strictly increasing if and only if $\theta(x) = \infty$.

As already discussed in the Introduction, it is convenient to write the energy dissipated by a crack opening (cf. Figure 1) as a function of the variation of the jump $V_u(t)$ defined in (2.3) (cf. Figure 2):

$$\int_\Gamma g(x, V_u(t)) \, \mathrm{d}\mathcal{H}^{n-1}.$$

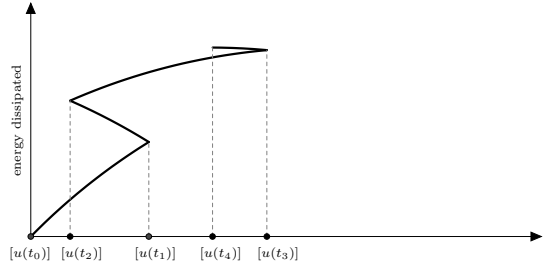

**Figure 1.** Energy dissipated by a jump $t \mapsto [u(t)]$ with a non-monotone history in a time interval $[t_0, t_4]$: $t \mapsto [u(t)]$ increases in $[t_0, t_1]$ and in $[t_2, t_3]$, whereas it decreases in $[t_1, t_2]$ and in $[t_3, t_4]$.

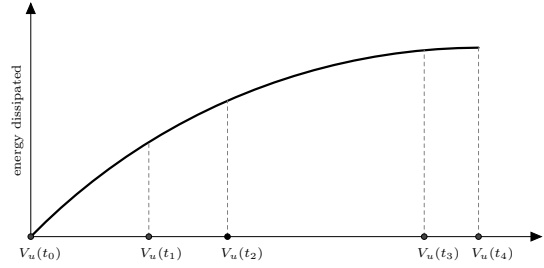

**Figure 2.** Energy dissipated as a function of the variation of the jumps $V_u(t)$ corresponding to a jump history as in Figure 1. Notice that the variation $V_u(t)$ is nondecreasing in time.

*Remark* 2.5. In [21] and [11], the variable used to describe the energy dissipated is the supremum of the jumps reached during the evolution. This is the main point where our cohesive model differs from those considered in [21] and [11].

**Definition of quasistatic evolution and strong formulation.** We are now in a position to give the definition of quasistatic evolution.

**Definition 2.6.** Let $w$, $u_0$, and $V_0$ be as in (2.2)–(2.5). Let $t \mapsto u(t)$ be a function defined on $[0, T]$ with values in $H^1(\Omega \setminus \Gamma)$ and let $V_u(t)$ be the variation of its jumps on $\Gamma$, defined in (2.3). We say that $t \mapsto u(t)$ is a *quasistatic evolution* with initial conditions $(u_0, V_0)$ and boundary datum $w$ if $u$ satisfies $u(0) = u_0$ and the following conditions:

(GS) *Global stability*: For every $t \in [0, T]$ we have $u(t) \in \mathscr{A}(w(t))$ and

$$\frac{1}{2} \int\limits_{\Omega \setminus \Gamma} |\nabla u(t)|^2 \, \mathrm{d}x + \int\limits_{\Gamma} g\big(x, V_u(t)\big) \, \mathrm{d}\mathcal{H}^{n-1} \leq \frac{1}{2} \int\limits_{\Omega \setminus \Gamma} |\nabla \widetilde{u}|^2 \, \mathrm{d}x + \int\limits_{\Gamma} g\big(x, \, V_u(t) + \big|[\widetilde{u}] - [u(t)]\big|\big) \, \mathrm{d}\mathcal{H}^{n-1} \,,$$

for every $\widetilde{u} \in \mathscr{A}(w(t))$.

(EB) *Energy-dissipation balance*: For every $t \in [0, T]$

$$\frac{1}{2} \int\limits_{\Omega \setminus \Gamma} |\nabla u(t)|^2 \, \mathrm{d}x + \int\limits_{\Gamma} g\big(x, V_u(t)\big) \, \mathrm{d}\mathcal{H}^{n-1}$$
$$= \frac{1}{2} \int\limits_{\Omega \setminus \Gamma} |\nabla u_0|^2 \, \mathrm{d}x + \int\limits_{\Gamma} g\big(x, V_0\big) \, \mathrm{d}\mathcal{H}^{n-1} + \int\limits_{0}^{t} \langle \nabla u(s), \nabla \dot{w}(s) \rangle_{L^2} \, \mathrm{d}s \,.$$

In order to give an insight into the strong formulation of the model studied in the paper, we state two results regarding necessary conditions satisfied by a quasistatic evolution. For simplicity, we derive these differential conditions under the assumption that $g(x, \cdot)$ is of class $\mathcal{C}^1$. We denote by $g'(x, \xi)$ the derivative of $g(x, \xi)$ with respect to $\xi$.

**Proposition 2.7.** *Assume that $g(x, \cdot)$ is of class $\mathcal{C}^1$ for $\mathcal{H}^{n-1}$-a.e. $x \in \Gamma$. Let $t \mapsto u(t)$ be a function defined on $[0, T]$ with values in $H^1(\Omega \setminus \Gamma)$ and satisfying* (GS). *Then for every $t \in [0, T]$ the following hold:*

(i) *The function $u(t)$ is a weak solution to the problem*

$$\begin{cases} \Delta u(t) = 0 & in \ \Omega \setminus \Gamma, \\ u(t) = w(t) & on \ \partial_D \Omega, \\ \partial_\nu u(t) = 0 & in \ H^{-\frac{1}{2}}(\partial_N \Omega). \end{cases}$$

(ii) *Let $u(t)^+ := u(t)|_{\Omega^+}$ and $u(t)^- := u(t)|_{\Omega^-}$. Then $\partial_\nu u(t)^+ = \partial_\nu u(t)^-$ in $H^{-\frac{1}{2}}(\Gamma)$.*

(iii) *Let $\partial_\nu u(t) := \partial_\nu u(t)^+ = \partial_\nu u(t)^-$. Then $\partial_\nu u(t) \in L^\infty(\Gamma)$ and*

$$|\partial_\nu u(t; x)| \le g'(x, V_u(t; x)) \quad for \ \mathcal{H}^{n-1}\text{-}a.e. \ x \in \Gamma. \tag{2.7}$$

To keep the presentation clear, the proof of Proposition 2.7 is given in Section 6.

Condition (iii) in Proposition 2.7 expresses the fact that the surface tension on $\Gamma$ due to the displacement is constrained to stay below a suitable threshold. This threshold decreases in time, since $g'(x, \cdot)$ is nonincreasing and $V_u(\cdot; x)$ is nondecreasing in time. However, this condition is static and is not enough to characterise an evolution.

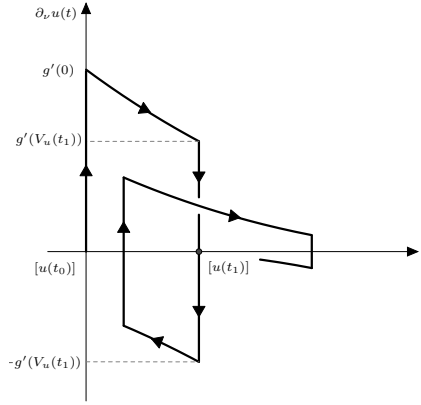

**Figure 3.** Crack opening versus surface tension corresponding to a jump history as in Figure 1.

Nonetheless, in the following proposition we employ the energy-dissipation balance to show that the evolution satisfies a flow rule: in the points where a crack opening grows, the surface tension actually must reach the maximal threshold. (See Figure 3 for a possible evolution of the surface tension.) The result is proved under regularity assumptions on the evolution $t \mapsto u(t)$. To make the statement concise, we denote by Sign the multifunction given by

$$\mathrm{Sign}(\xi) := \begin{cases} 1 & \text{if } \xi > 0, \\ [-1, 1] & \text{if } \xi = 0, \\ -1 & \text{if } \xi < 0. \end{cases}$$

**Proposition 2.8.** *Assume that $g(x, \cdot)$ is of class $\mathcal{C}^1$ for $\mathcal{H}^{n-1}$-a.e. $x \in \Gamma$. Let $t \mapsto u(t)$ be a quasistatic evolution in the sense of Definition 2.6 and assume that $u \in AC([0, T]; H^1(\Omega \setminus \Gamma))$. Then*

$$\partial_\nu u(t; x) \in g'(x, V_u(t; x)) \, \mathrm{Sign}([\dot{u}(t; x)]) \quad for \ \mathcal{H}^{n-1}\text{-}a.e. \ x \in \Gamma \ and \ a.e. \ t \in [0, T],$$

*where $[\dot{u}(t)]$ is the derivative in time of $[u(t)]$ with respect to the strong topology in $L^2(\Gamma)$.*

Proposition 2.8 is proved in Section 6.

**Statement of the main result.** We now introduce the tools needed to state our main result, which concern the existence of a quasistatic evolution and the approximation by means of discrete-time evolutions.

As usual in the proof of existence of quasistatic evolutions for rate-independent systems, we construct discrete-time evolutions by solving incremental minimum problems. For every $k \in \mathbb{N}$, let us consider a subdivision of the time interval $[0, T]$ given by $k+1$ nodes

$$0 = t_k^0 < t_k^1 < \cdots < t_k^{k-1} < t_k^k = T, \quad \lim_{k \to \infty} \max_{1 \leq i \leq k} |t_k^i - t_k^{i-1}| = 0,$$

and let us define $w_k^i := w(t_k^i)$.

We assume that the initial condition $(u_0, V_0)$ is globally stable, namely

$$\frac{1}{2} \int_{\Omega \setminus \Gamma} |\nabla u_0|^2 \, \mathrm{d}x + \int_\Gamma g(x, V_0) \, \mathrm{d}\mathcal{H}^{n-1} \leq \frac{1}{2} \int_{\Omega \setminus \Gamma} |\nabla \widetilde{u}|^2 \, \mathrm{d}x + \int_\Gamma g\left(x, V_0 + |[\widetilde{u}] - [u_0]|\right) \mathrm{d}\mathcal{H}^{n-1}, \quad (2.8)$$

for every $\widetilde{u} \in \mathscr{A}(w(0))$.

As the first step of the incremental process, we set $u_k^0 := u_0$ and $V_k^0 := V_0$. Let $i \in \{1, \ldots, k\}$ and assume that we know $u_k^h$ and $V_k^h$ for $h = 0, \ldots, i-1$. Then we define $u_k^i$ as a solution to the problem

$$\min_u \left\{ \frac{1}{2} \int_{\Omega \setminus \Gamma} |\nabla u|^2 \, \mathrm{d}x + \int_\Gamma g\left(x, V_k^{i-1} + |[u] - [u_k^{i-1}]|\right) \mathrm{d}\mathcal{H}^{n-1} \ : \ u \in \mathscr{A}(w_k^i) \right\}, \quad (2.9)$$

and we set

$$V_k^i := V_k^{i-1} + |[u_k^i] - [u_k^{i-1}]| = V_0 + \sum_{j=1}^i |[u_k^j] - [u_k^{j-1}]|. \quad (2.10)$$

The existence of a solution to (2.9) is obtained by employing the direct method of the Calculus of Variations.

The discrete-time evolutions are then defined as piecewise constant interpolations of the solutions to the incremental problems. Namely, we set

$$u_k(t) := u_k^i, \quad V_k(t) := V_k^i, \quad w_k(t) := w_k^i \quad \text{for } t_k^i \leq t < t_k^{i+1} \quad (2.11)$$

and $u_k(T) := u_k^k$, $V_k(T) := V_k^k$, $w_k(T) := w(T)$.

Passing to the limit as $k \to \infty$, we prove that $u_k$ converges to a quasistatic evolution $u$. A major point of our result is that the convergence holds for a subsequence independent of $t$. We also provide a convergence result for the variations of the jumps. Specifically, the truncated functions $V_k(t) \wedge \theta$ converge to $V_u(t) \wedge \theta$, where $\theta$ is as in (2.6), and $\wedge$ denotes the minimum between two functions. We remark that when $V_u(t)$ overcomes the threshold $\theta(x)$, we have no control on $V_u(t; x)$, which may increase without further dissipation of energy. Moreover, we obtain that $t \mapsto u(t)$ and $t \mapsto V_u(t)$ are continuous (in a suitable sense), except for countably many times. These results are stated in the following theorem, whose proof is given in Section 6.

**Theorem 2.9** (Existence and approximation of quasistatic evolutions)**.** *Assume that $g$ satisfies* $(g1)-(g4)$. *Let $w$, $u_0$, and $V_0$ be as in* (2.2)-(2.5) *and assume that $(u_0, V_0)$ is globally stable in the sense of* (2.8). *Consider the piecewise constant evolutions $t \mapsto u_k(t)$ and the piecewise constant variations $t \mapsto V_k(t)$ defined in* (2.11). *Then there exist a subsequence (independent of $t$ and not relabelled) and a quasistatic evolution $t \mapsto u(t)$ with initial conditions $(u_0, V_0)$ and boundary datum $w$ such that, for every $t \in [0, T]$,*

$$u_k(t) \to u(t) \qquad \text{strongly in } H^1(\Omega \setminus \Gamma), \quad (2.12)$$

$$V_k(t) \wedge \theta \to V_u(t) \wedge \theta \quad \text{in measure}, \quad (2.13)$$

*where $V_u(t)$ is the function defined in* (2.3) *and $\theta$ is given in* (2.6).

*Moreover, there exists a set $E \subset [0, T]$, at most countable, such that, for every $t \in [0, T] \setminus E$ and every $s \to t$,*

$$u(s) \to u(t) \qquad \text{strongly in } H^1(\Omega \setminus \Gamma). \quad (2.14)$$

$$V_u(s) \wedge \theta \to V_u(t) \wedge \theta \quad \text{in measure}. \quad (2.15)$$

We underline that, if $\theta(x)$ is finite and $V_u(t; x) \geq \theta(x)$, the material is completely broken at $x$. Therefore $V_u(t) \wedge \theta$, appearing in the theorem above, is the relevant state variable for the system.

*Remark* 2.10. If $\theta \in L^\infty(\Gamma)$, then the convergence in (2.13) and (2.15) is also strong in $L^p(\Gamma)$ for every $p \in [1, \infty)$. In contrast, if $\theta \equiv \infty$ (that is $g(x, \cdot)$ is strictly increasing for $\mathcal{H}^{n-1}$-a.e. $x \in \Gamma$), then $V_k(t) \to V_u(t)$ in measure as $k \to \infty$ and $V_k(s) \to V_u(t)$ in measure as $s \to t$.

**Guidelines for the proof of the main result.** The main difficulty in the passage to the continuous-time limit as $k \to \infty$ is that we lack of controls on $V_k(t)$. In fact, by (2.9), we can only infer that $\int_\Gamma g(x, V_k(t)) \, d\mathcal{H}^{n-1}$ is uniformly bounded, but this gives no information on $V_k(t)$, since $g$ is bounded. For this reason we resort to a weaker notion of quasistatic evolution, where the variation of jumps on $\Gamma$ is replaced by a Young measure. Notwithstanding, after establishing the properties of such an evolution, we are able to show that the Young measure found in the limit is concentrated on a function. Eventually, we obtain a quasistatic evolution in the sense of Definition 2.6. We describe here the strategy followed to prove Theorem 2.9.

Following the scheme of the proof of existence of energetic solutions to rate-independent systems [27], the starting point of our analysis is to obtain a global stability and an energy-dissipation inequality for the discrete-time evolutions $t \mapsto u_k(t)$ (Proposition 4.1). As usual, the energy-dissipation inequality provides a priori bounds in $H^1(\Omega \setminus \Gamma)$ for the functions $u_k(t)$, independently of $k$ and $t$. In order to study the limit of the functions $V_k(t)$, it is convenient to introduce the Young measures concentrated on the graph of $V_k(t)$, namely

$$\nu_k(t) := \delta_{V_k(t)} \in \mathcal{Y}(\Gamma; [0, \infty]) \quad \text{for every } t \in [0, T] \, . \tag{2.16}$$

We refer to Section 3 for the notation and the basic properties of Young measures. Since the functions $V_k(t)$ are nondecreasing with respect to $t$, we can apply a Helly-type selection principle (proved in [11]) to infer that the Young measures $\nu_k(t)$ converge narrowly to a Young measure $\nu(t) \in \mathcal{Y}(\Gamma; [0, \infty])$ on a subsequence independent of $t$. Thanks to the a priori bounds on $u_k(t)$, it is possible to extract a subsequence $k_j(t)$ (depending on $t$) such that $u_{k_j(t)}(t)$ converges to $u(t)$ weakly in $H^1(\Omega \setminus \Gamma)$. These convergences allow us to pass to the limit in the global stability of the discrete-time evolutions (Proposition 4.4), and thus to deduce that $t \mapsto (u(t), \nu(t))$ satisfies a suitable notion of global stability (condition (GSY) in Definition 5.1).

Afterwards, we show that the evolution $t \mapsto (u(t), \nu(t))$ satisfies an energy-dissipation balance (condition (EBY) in Definition 5.1). One inequality in this balance is a consequence of the energy-dissipation inequality of the discrete-time evolutions $t \mapsto u_k(t)$. On the contrary, the proof of the opposite inequality requires a thorough analysis. The main reason is that the Helly Selection Principle adopted before does not give any information about the relation between the Young measure $\nu(t)$ and $V_u(t)$. This relation is though encoded in a property satisfied by $t \mapsto \nu(t)$ (the *irreversibility* condition (IRY) in Definition 5.1), that we derive from the analogous condition $(IRY)_k$ for the approximating Young measures $t \mapsto \nu_k(t)$. This property relates $\nu(t)$ to $[u(t)]$ and allows us to conclude the proof of the other inequality in the energy-dissipation balance by employing the global stability.

In addition, we prove that $u_k(t)$ actually converges to $u(t)$ strongly in $H^1(\Omega \setminus \Gamma)$ on a subsequence independent of $t$. This convergence result is proved in Section 5 by showing that the jump $\gamma(t) := [u(t)]$ is determined *de facto* independently of $t$ (cf. equation (5.9)). Indeed this implies that the function $u(t)$ is the unique solution of a minimum problem among functions with a prescribed jump $\gamma(t)$ (Proposition 5.6). With similar arguments, we prove that $t \mapsto u(t)$ is continuous in $t$ except for a countable set $E \subset [0, T]$.

Finally, in Section 6 we prove that $u$ is actually a quasistatic evolution in the sense of Definition 2.6. Notice that for this step we need the assumption on the concavity of $g(x, \cdot)$. Moreover, this allows us to prove that the Young measure $\nu(t)$ (suitably truncated with $\theta$) is concentrated on the function $V_u(t)$. As a consequence of this fact, we are able to deduce also the convergences in (2.13) and (2.15) in Theorem 2.9.

## 3. Preliminary results about Young measures

**Probability measures.** Let $\Xi$ be a metric space. We denote by $\mathcal{M}_b^+(\Xi)$ the set of positive bounded measures, and by $\mathcal{P}(\Xi)$ the set of probability measures on $\Xi$. The space $\mathcal{M}_b^+(\Xi)$ can be put in duality with the space of bounded continuous functions $\mathcal{C}_b(\Xi)$ by defining

$$\langle f, \mu \rangle := \int_\Xi f(\xi) \mu(d\xi) = \int_\Xi f(\xi) \, d\mu(\xi) \, , \tag{3.1}$$

for every $\mu \in \mathcal{M}_b^+(\Xi)$ and $f \in \mathcal{C}_b(\Xi)$.

If $\Xi$ is a separable metric space and $\mu \in \mathcal{M}_b^+(\Xi)$, the support of $\mu$ is the smallest closed subset of $\Xi$ where the measure $\mu$ is concentrated, i.e.,

$$\operatorname{supp}(\mu) := \bigcap_{\substack{C \text{ closed} \\ \mu(\Xi \setminus C) = 0}} C \,.$$

Let $\Xi_1$ and $\Xi_2$ be two metric spaces, let $\varphi \colon \Xi_1 \to \Xi_2$ be a Borel map, and let $\mu \in \mathcal{M}_b^+(\Xi_1)$. The *push-forward* of $\mu$ through the map $\varphi$ is the measure $\varphi_\# \mu \in \mathcal{M}_b^+(\Xi_2)$ defined by $\varphi_\# \mu(A) := \mu(\varphi^{-1}(A))$ for every $A \in \mathcal{B}(\Xi_2)$.

We will later deal with measures in the space $\mathcal{M}_b^+([-\infty, \infty])$, where $[-\infty, \infty]$ is endowed with the metric induced by an increasing homeomorphism

$$\phi \colon [-\infty, \infty] \to [-1, 1] \,, \tag{3.2}$$

e.g. $\phi(\xi) := \frac{2}{\pi} \arctan(\xi)$. Measures in $\mathcal{M}_b^+([-\infty, \infty])$ are in duality with bounded continuous functions $f \in \mathcal{C}_b([-\infty, \infty])$, i.e., continuous functions with a finite limit at $\pm\infty$.

We also recall that for every probability measure $\mu \in \mathcal{P}([-\infty, \infty])$ we can define the *cumulative distribution function* $F_\mu \colon [-\infty, \infty] \to [0, 1]$ by

$$F_\mu(\xi) := \mu([-\infty, \xi]) \quad \text{for every } \xi \in [-\infty, \infty] \,. \tag{3.3}$$

By the right continuity of $F_\mu$, it is possible to define its *pseudo-inverse* $F_\mu^{[-1]} \colon [0, 1] \to [-\infty, \infty]$ by

$$F_\mu^{[-1]}(m) := \min\{\xi \in \mathbb{R} : F_\mu(\xi) \geq m\} \,. \tag{3.4}$$

**Young measures.** For an introduction to the general theory of Young measures we refer, e.g., to [28]. Here we recall some basic notions and properties. Let us fix a metric space $\Xi$.

**Definition 3.1.** The collection of *Young measures* on $\Gamma \times \Xi$ with respect to the measure $\mathcal{H}^{n-1}$ is the set

$$\mathcal{Y}(\Gamma; \Xi) := \{\nu \in \mathcal{M}_b^+(\Gamma \times \Xi) : \pi_\#^\Gamma \nu = \mathcal{H}^{n-1} \llcorner \Gamma\} \,,$$

where $\pi^\Gamma \colon \Gamma \times \Xi \to \Gamma$ is the projection on $\Gamma$.

*Remark* 3.2. We recall that a family $(\nu^x)_{x \in \Gamma}$ of probability measures $\nu^x \in \mathcal{P}(\Xi)$ parametrised on $\Gamma$ is said to be *measurable* if the function $x \mapsto \nu^x(A)$ is $\mathcal{H}^{n-1}$-measurable for every $A \in \mathcal{B}(\Xi)$. By the Disintegration Theorem (see [5, Theorem 2.28]), it is always possible to associate a measurable family of probability measures $(\nu^x)_{x \in \Gamma}$ with a Young measure $\nu \in \mathcal{Y}(\Gamma; X)$ in such a way that

$$\int_{\Gamma \times \Xi} f(x, \xi) \, \mathrm{d}\nu = \int_\Gamma \int_\Xi f(x, \xi) \nu^x(\mathrm{d}\xi) \, \mathrm{d}\mathcal{H}^{n-1} \quad \text{for every } f \in L^1_\nu(\Gamma \times X) \,. \tag{3.5}$$

Moreover, the family $(\nu^x)_{x \in \Gamma}$ is unique up to $\mathcal{H}^{n-1}$-negligible sets, i.e., if $(\widetilde{\nu}^x)_{x \in \Gamma}$ is any other measurable family of probability functions satisfying (3.5), then $\widetilde{\nu}^x = \nu^x$ for $\mathcal{H}^{n-1}$-a.e. $x \in \Gamma$.

If $\nu = (\nu^x)_{x \in \Gamma} \in \mathcal{Y}(\Gamma; \Xi)$, for every $f \in \mathcal{C}_b(\Gamma \times \Xi)$ the duality between $\nu$ and $f$ reads

$$\int_{\Gamma \times \Xi} f(x, \xi) \, \mathrm{d}\nu(x, \xi) = \int_\Gamma \int_\Xi f(x, \xi) \nu^x(\mathrm{d}\xi) \, \mathrm{d}\mathcal{H}^{n-1} = \int_\Gamma \langle f(x, \cdot), \nu^x \rangle \, \mathrm{d}\mathcal{H}^{n-1} \,.$$

*Example* 3.3. The simplest example of a Young measure is obtained by fixing a measurable function $v \colon \Gamma \to \Xi$ and by considering the Young measure *concentrated* on the graph of the function $v$, identified by the measurable family of probability measures $\delta_v := (\delta_{v(x)})_{x \in \Gamma}$.

We will consider the space $\mathcal{Y}(\Gamma; \Xi)$ endowed with the *narrow* topology.

**Definition 3.4.** We say that $\nu_j$ converges *narrowly* to $\nu$ (and denote $\nu_j \rightharpoonup \nu$) if and only if

$$\int_\Gamma \langle f(x, \cdot), \nu_j^x \rangle \, \mathrm{d}\mathcal{H}^{n-1} \to \int_\Gamma \langle f(x, \cdot), \nu^x \rangle \, \mathrm{d}\mathcal{H}^{n-1}, \tag{3.6}$$

for every $f \in \mathcal{C}_b(\Gamma \times \Xi)$.

*Remark* 3.5. If $\Xi$ is a compact metric space, by [28, Theorem 2] the convergence in (3.6) also holds for every *Carathéodory integrand* $f$, i.e., a measurable function such that $f(x, \cdot) \in \mathcal{C}_b(\Xi)$ for $\mathcal{H}^{n-1}$-a.e. $x \in \Gamma$ and such that $x \mapsto \|f(x, \cdot)\|_\infty$ belongs to $L^1(\Gamma)$.

The narrow convergence for concentrated Young measures is characterised in the following proposition. For the proof, we refer to [28, Proposition 6].

**Proposition 3.6.** *Assume that $\Xi$ is a compact metric space. Let $v_j, v \colon \Gamma \to \Xi$ be measurable functions. Then $\delta_{v_j} \to \delta_v$ if and only if $v_j \to v$ in measure.*

*Remark* 3.7. In the case where $\Xi$ is $[-\infty, \infty]$ endowed with the metric induced by $\phi$ in (3.2), then $v_j \to v$ in measure if and only if $\mathcal{H}^{n-1}(\{|\phi(v_j) - \phi(v)| \geq \varepsilon\}) \to 0$ for every $\varepsilon > 0$.

The following compactness result holds (cf. [28, Theorem 2]).

**Theorem 3.8.** *Assume that $\Xi$ is a compact metric space. Then $\mathcal{Y}(\Gamma; \Xi)$, endowed with the narrow topology, is sequentially compact.*

*Remark* 3.9. The assumption on the compactness of the space $\Xi$ is crucial to guarantee the compactness of $\mathcal{Y}(\Gamma; \Xi)$ with respect to the narrow convergence. For instance, if $\Xi = \mathbb{R}$, it may happen that a sequence $\nu_j \in \mathcal{Y}(\Gamma; \mathbb{R})$ has some mass escaping to infinity.

We will later need to infer the compactness of sequences $\nu_j \in \mathcal{Y}(\Gamma; \mathbb{R})$ with no tightness assumptions. Thus, we will consider a compactification of $\mathbb{R}$, i.e., we will regard $\nu_j$ as Young measures in $\mathcal{Y}(\Gamma; [-\infty, \infty])$. In this way, we can conclude that a subsequence of $\nu_j \in \mathcal{Y}(\Gamma; \mathbb{R})$ (not relabelled) converges narrowly to a Young measure $\nu \in \mathcal{Y}(\Gamma; [-\infty, \infty])$.

To deal with these Young measures, it is convenient to introduce the map

$$\Phi \colon \Gamma \times [-\infty, \infty] \to \Gamma \times [-1, 1], \quad \Phi(x, \xi) := (x, \phi(x)), \tag{3.7}$$

where $\phi$ is the homeomorphism defined in (3.2). In this way, for every $\nu \in \mathcal{Y}(\Gamma; [-\infty, \infty])$ we have $\Phi_\# \nu \in \mathcal{Y}(\Gamma; [-1, 1])$. The elements of $\mathcal{Y}(\Gamma; [-\infty, \infty])$ are in duality with functions $f \in \mathcal{C}_b(\Gamma \times [-\infty, \infty])$, i.e., such that $f \circ \Phi^{-1} \in \mathcal{C}_b(\Gamma \times [-1, 1])$.

**Translation.** We now recall how to shift real-valued Young measures. For every measurable function $\gamma \colon \Gamma \to \mathbb{R}$ we define the translation map $\mathcal{S}^\gamma \colon \Gamma \times [-\infty, \infty] \to \Gamma \times [-\infty, \infty]$ by $\mathcal{S}^\gamma(x, \xi) := (x, \xi + \gamma(x))$, with the usual convention that $a \pm \infty = \pm\infty$ for every $a \in \mathbb{R}$. For every $\nu \in \mathcal{Y}(\Gamma; [-\infty, \infty])$ we set

$$\nu \oplus \gamma := \mathcal{S}^\gamma_\# \nu \in \mathcal{Y}(\Gamma; [-\infty, \infty]), \tag{3.8}$$

$$\nu \ominus \gamma := \mathcal{S}^{(-\gamma)}_\# \nu \in \mathcal{Y}(\Gamma; [-\infty, \infty]). \tag{3.9}$$

*Remark* 3.10. Let $\nu_j, \nu \in \mathcal{Y}(\Gamma; [-\infty, \infty])$ be such that $\nu_j \rightharpoonup \nu$ and let $\gamma \colon \Gamma \to \mathbb{R}$ be a measurable function. By Remark 3.5 we have $\nu_j \oplus \gamma \rightharpoonup \nu \oplus \gamma$.

Moreover, if $\gamma, \gamma_j \colon \Gamma \to \mathbb{R}$ are such that $\gamma_j \to \gamma$ in measure, then it is easy to see that $\nu_j \oplus \gamma_j \rightharpoonup \nu \oplus \gamma$.

**Truncation.** We now introduce the notion of truncation of Young measures. This will be employed in Section 6. Given a Young measure $\nu \in \mathcal{Y}(\Gamma; [-\infty, \infty])$ and a measurable function $\theta \colon \Gamma \to [-\infty, \infty]$, we consider the map $\mathcal{T}^\theta \colon \Gamma \times [-\infty, \infty] \to \Gamma \times [-\infty, \infty]$ given by

$$\mathcal{T}^\theta(x, \xi) := (x, \xi \wedge \theta(x)) \tag{3.10}$$

and we say that $\mathcal{T}^\theta_\# \nu$ is the *truncation* of $\nu$ by $\theta$.

*Remark* 3.11. In this case, the cumulative distribution function of the measure $(\mathcal{T}^\theta_\# \nu)^x$ is given by

$$F_{(\mathcal{T}^\theta_\# \nu)^x}(\xi) = \begin{cases} F_{\nu^x}(\xi) & \text{if } \xi < \theta(x), \\ 1 & \text{if } \xi \geq \theta(x), \end{cases}$$

for $\mathcal{H}^{n-1}$-a.e. $x \in \Gamma$. Moreover, if $\nu_j \rightharpoonup \nu$ in $\mathcal{Y}(\Gamma; [-\infty, \infty])$, then by Remark 3.5 we have $\mathcal{T}^\theta_\# \nu_j \rightharpoonup \mathcal{T}^\theta_\# \nu$ in $\mathcal{Y}(\Gamma; [-\infty, \infty])$.

**Partial order.** Following [11, Definition 3.10], we introduce a partial order in the space of Young measures on $\Gamma \times \mathbb{R}$. We recall here the definition of this order and its main properties.

**Definition 3.12.** Let $\nu_1 = (\nu_1^x)_{x \in \Gamma}$, $\nu_2 = (\nu_2^x)_{x \in \Gamma} \in \mathcal{Y}(\Gamma; \mathbb{R})$. We say that $\nu_1 \preceq \nu_2$ if one of the following equivalent conditions is satisfied:

(i) for every Carathéodory integrand $f \colon \Gamma \times \mathbb{R} \to \mathbb{R}$ nondecreasing with respect to the second variable we have

$$\int_\Gamma \langle f(x, \cdot), \nu_1^x \rangle \, \mathrm{d}\mathcal{H}^{n-1} \leq \int_\Gamma \langle f(x, \cdot), \nu_2^x \rangle \, \mathrm{d}\mathcal{H}^{n-1};$$

(ii) $F_{\nu_1^x}(\xi) \geq F_{\nu_2^x}(\xi)$ for $\mathcal{H}^{n-1}$-a.e. $x \in \Gamma$ and for every $\xi \in \mathbb{R}$.

*Remark* 3.13. If $\nu_1$ and $\nu_2$ are concentrated on some measurable functions $\gamma_1$ and $\gamma_2$, respectively, then

$$\nu_1 \preceq \nu_2 \qquad \text{if and only if} \qquad \gamma_1(x) \leq \gamma_2(x) \text{ for } \mathcal{H}^{n-1}\text{-a.e. } x \in \Gamma.$$

The partial order $\preceq$ is naturally extended to Young measures $\mathcal{Y}(\Gamma; [-\infty, \infty])$ by employing the homeomorphism $\Phi$ defined in (3.7). Namely, for every $\nu_1$, $\nu_2 \in \mathcal{Y}(\Gamma; [-\infty, \infty])$ we have $\nu_1 \preceq \nu_2$ if and only if $\Phi_\# \nu_1 \preceq \Phi_\# \nu_2$.

In the following we recall the definition of supremum of a family of Young measures. (See [11, Proposition 3.16] for the existence of such a Young measure.)

**Definition 3.14.** Let $(\nu_i)_{i \in I}$ be a family of Young measures in $\mathcal{Y}(\Gamma; [-\infty, \infty])$. We say that $\overline{\nu} \in \mathcal{Y}(\Gamma; [-\infty, \infty])$ is the *supremum* over $i \in I$ of the family $(\nu_i)_{i \in I}$, and we write

$$\overline{\nu} = \sup_{i \in I} \nu_i,$$

if the following two conditions hold:

(i) $\overline{\nu} \succeq \nu_i$ for every $i \in I$;
(ii) if $\nu \in \mathcal{Y}(\Gamma; [-\infty, \infty])$ such that $\nu \succeq \nu_i$ for every $i \in I$, then $\nu \succeq \overline{\nu}$.

*Remark* 3.15. In the case where $\nu_i$ are concentrated on measurable functions $v_i \colon \Gamma \to [-\infty, \infty]$, $i \in I$, we have

$$\sup_{i \in I} \delta_{v_i} = \delta_{\overline{v}},$$

where $\overline{v} = \operatorname*{ess\,sup}_{i \in I} v_i$ (cf. [11, Remark 3.17]).

*Remark* 3.16. If a map $t \mapsto \nu(t)$ from $[0, T]$ to $\mathcal{Y}(\Gamma; [-\infty, \infty])$ is nondecreasing with respect to $\preceq$, then there exists a countable set $E \subset [0, T]$ such that $t \mapsto \nu(t)$ is continuous in $[0, T] \setminus E$. The proof of this fact is an easy consequence of [11, Lemma 3.19].

We conclude this section by recalling the Helly Selection Principle for Young measures [11, Theorem 3.20], a key tool for the proof of our result. Notice that [11, Theorem 3.20] is stated for Young measures with values in $\mathbb{R}$ instead of $[-\infty, \infty]$.

**Theorem 3.17.** *Let $t \mapsto \nu_k(t)$ be a sequence of maps from $[0, T]$ to $\mathcal{Y}(\Gamma; [-\infty, \infty])$ that are nondecreasing with respect to $\preceq$. Then there exists a subsequence $\nu_{k_j}$, independent of $t$, and a nondecreasing map $t \mapsto \nu(t)$ from $[0, T]$ to $\mathcal{Y}(\Gamma; [-\infty, \infty])$ such that $\nu_{k_j}(t) \rightharpoonup \nu(t)$, as $j \to \infty$, for every $t \in [0, T]$.*

*Proof.* The result follows from a straightforward application of [11, Theorem 3.20] to the sequence of nondecreasing maps $\Phi_\# \nu_k(t) \in \mathcal{Y}(\Gamma; [-1, 1])$, where $\Phi$ is the homeomorphism $\Phi$ defined in (3.7). $\square$

## 4. DISCRETE-TIME EVOLUTIONS

We study here the discrete-time evolutions already introduced in Section 2.

Let $u_k(t)$, $V_k(t)$, and $w_k(t)$ be the piecewise constant interpolations given in (2.11). Let $\nu_k(t) \in \mathcal{Y}(\Gamma; [0, \infty])$ be the Young measures concentrated on $V_k(t)$ defined in (2.16). In the following proposition we state the main properties satisfied by such approximate evolutions and we provide a priori bounds for $u_k(t)$.

**Proposition 4.1.** *The discrete evolutions $t \mapsto u_k(t)$ defined in (2.11) satisfy the following conditions:*

$(\mathrm{GS})_k$  *Global stability: For every $t \in [0, T]$ we have $u_k(t) \in \mathscr{A}(w_k(t))$ and*

$$\frac{1}{2}\int\limits_{\Omega\setminus\Gamma} |\nabla u_k(t)|^2 \,\mathrm{d}x + \int\limits_{\Gamma} g(x, V_k(t)) \,\mathrm{d}\mathcal{H}^{n-1} \leq \frac{1}{2}\int\limits_{\Omega\setminus\Gamma} |\nabla \widetilde{u}|^2 \,\mathrm{d}x + \int\limits_{\Gamma} g\big(x,\, V_k(t) + \big|[\widetilde{u}] - [u_k(t)]\big|\big) \,\mathrm{d}\mathcal{H}^{n-1},$$

*for every $\widetilde{u} \in \mathscr{A}(w_k(t))$.*

$(\mathrm{EI})_k$  *Energy-dissipation inequality: There exists a sequence $\eta_k$ with $\eta_k \to 0$ as $k \to \infty$ such that for every $t \in [0, T]$ we have*

$$\frac{1}{2}\int\limits_{\Omega\setminus\Gamma} |\nabla u_k(t)|^2 \,\mathrm{d}x + \int\limits_{\Gamma} g(x, V_k(t)) \,\mathrm{d}\mathcal{H}^{n-1}$$

$$\leq \frac{1}{2}\int\limits_{\Omega\setminus\Gamma} |\nabla u_0|^2 \,\mathrm{d}x + \int\limits_{\Gamma} g(x, V_0) \,\mathrm{d}\mathcal{H}^{n-1} + \int\limits_0^{t_k^i} \langle \nabla u_k(s), \nabla \dot{w}(s) \rangle_{L^2} \,\mathrm{d}s + \eta_k \,,$$

*where $i \in \{0, \dots, k\}$ is the largest integer such that $t_k^i \leq t$.*

*Moreover, there exists a constant $C > 0$ independent of $k$ and $t$ such that*

$$\|u_k(t)\|_{H^1(\Omega\setminus\Gamma)} \leq C \quad \text{for every } k \in \mathbb{N} \text{ and } t \in [0, T]\,. \tag{4.1}$$

*Proof.* In order to prove the global stability $(\mathrm{GS})_k$, we notice that if $i$ is the largest integer such that $t_k^i \leq t$, then by (2.10) we get that

$$V_k(t) + \big|[\widetilde{u}] - [u_k(t)]\big| = V_k^i + \big|[\widetilde{u}] - [u_k^i]\big| = V_k^{i-1} + \big|[u_k^i] - [u_k^{i-1}]\big| + \big|[\widetilde{u}] - [u_k^i]\big|$$
$$\geq V_k^{i-1} + \big|[\widetilde{u}] - [u_k^{i-1}]\big|\,.$$

Then we infer $(\mathrm{GS})_k$ by the fact that $u_k(t) = u_k^i$ is a solution to (2.9) and by the monotonicity of $g(x, \cdot)$.

Let us prove the energy-dissipation inequality $(\mathrm{EI})_k$. Let us fix $t \in [0, T]$, $k \in \mathbb{N}$, and $i \in \{1, \dots, k\}$ as in the statement (the case $i = 0$ being trivial). For $1 \leq h \leq i$, the function $u_k^{h-1} - w_k^{h-1} + w_k^h$ is an admissible competitor for the minimum problem (2.9) solved by $u_k^h$. Hence

$$\frac{1}{2}\int\limits_{\Omega\setminus\Gamma} |\nabla u_k^h|^2 \,\mathrm{d}x + \int\limits_{\Gamma} g(x, V_k^h) \,\mathrm{d}\mathcal{H}^{n-1} \leq \frac{1}{2}\int\limits_{\Omega\setminus\Gamma} |\nabla u_k^{h-1}|^2 \,\mathrm{d}x + \int\limits_{\Gamma} g(x, V_k^{h-1}) \,\mathrm{d}\mathcal{H}^{n-1}$$

$$+ \int\limits_{\Omega\setminus\Gamma} \nabla u_k^{h-1} \cdot (\nabla w_k^h - \nabla w_k^{h-1}) \,\mathrm{d}x + \frac{1}{2}\int\limits_{\Omega\setminus\Gamma} |\nabla w_k^h - \nabla w_k^{h-1}|^2 \,\mathrm{d}x$$

$$\leq \frac{1}{2}\int\limits_{\Omega\setminus\Gamma} |\nabla u_k^{h-1}|^2 \,\mathrm{d}x + \int\limits_{\Gamma} g(x, V_k^{h-1}) \,\mathrm{d}\mathcal{H}^{n-1}$$

$$+ \int\limits_{t_k^{h-1}}^{t_k^h} \langle \nabla u_k(s), \nabla \dot{w}(s) \rangle_{L^2} \,\mathrm{d}s + \frac{1}{2}\bigg( \int\limits_{t_k^{h-1}}^{t_k^h} \|\nabla \dot{w}(s)\|_{L^2} \,\mathrm{d}s \bigg)^2,$$

$$\tag{4.2}$$

where we used our assumption (2.2) on $w$ to deduce that

$$\nabla w_k^h - \nabla w_k^{h-1} = \int\limits_{t_k^{h-1}}^{t_k^h} \nabla \dot{w}(s) \,\mathrm{d}s\,,$$

as a Bochner integral in $L^2$. Summing up the inequalities given by (4.2) for $h = 1, \dots, i$, we get $(\mathrm{EI})_k$ with

$$\eta_k := \frac{1}{2} \bigg( \max_{1 \leq h \leq k} \int\limits_{t_k^{h-1}}^{t_k^h} \|\nabla \dot{w}(s)\|_{L^2} \,\mathrm{d}s \bigg) \bigg( \int\limits_0^T \|\nabla \dot{w}(s)\|_{L^2} \,\mathrm{d}s \bigg)\,.$$

In particular, from $(\mathrm{EI})_k$ we readily deduce that there exists a constant $C > 0$ independent of $k$ and $t$ such that $\|\nabla u_k(t)\|_{L^2} \leq C$. Then, by the Poincaré inequality, we get (4.1) (up to changing the name of the constant). $\qquad \square$

*Remark* 4.2. It is convenient to express the properties satisfied by $u_k(t)$ also in terms of the Young measures $\nu_k(t) \in \mathcal{Y}(\Gamma; [0, \infty])$ defined in (2.16). In Section 5, we will pass to the limit in these conditions.

$\mathrm{(IRY)}_k$  *Irreversibility*: $\nu_k(t) \succeq \nu_k(s) \oplus \big|[u_k(t)] - [u_k(s)]\big|$ for every $s, t \in [0, T]$ with $s \leq t$.

$\mathrm{(GSY)}_k$  *Global stability*: For every $t \in [0, T]$ we have $u_k(t) \in \mathscr{A}(w_k(t))$ and

$$\frac{1}{2} \int\limits_{\Omega \backslash \Gamma} |\nabla u_k(t)|^2 \, \mathrm{d}x + \int\limits_{\Gamma} \langle g(x, \cdot), \nu_k^x(t) \rangle \, \mathrm{d}\mathcal{H}^{n-1} \leq \frac{1}{2} \int\limits_{\Omega \backslash \Gamma} |\nabla \widetilde{u}|^2 \, \mathrm{d}x + \int\limits_{\Gamma} \langle g(x, \cdot), \widetilde{\nu}_k^x \rangle \, \mathrm{d}\mathcal{H}^{n-1},$$

for every $\widetilde{u} \in \mathscr{A}(w_k(t))$, where $\widetilde{\nu}_k := \nu_k(t) \oplus \big|[\widetilde{u}] - [u_k(t)]\big| \in \mathcal{Y}(\Gamma; [0, \infty])$.

$\mathrm{(EIY)}_k$  *Energy-dissipation inequality*: For every $t \in [0, T]$

$$\frac{1}{2} \int\limits_{\Omega \backslash \Gamma} |\nabla u_k(t)|^2 \, \mathrm{d}x + \int\limits_{\Gamma} \langle g(x, \cdot), \nu_k^x(t) \rangle \, \mathrm{d}\mathcal{H}^{n-1}$$
$$\leq \frac{1}{2} \int\limits_{\Omega \backslash \Gamma} |\nabla u_0|^2 \, \mathrm{d}x + \int\limits_{\Gamma} g(x, V_0) \, \mathrm{d}\mathcal{H}^{n-1} + \int\limits_0^{t_k^i} \langle \nabla u_k(s), \nabla \dot{w}(s) \rangle_{L^2} \, \mathrm{d}s + \eta_k \,,$$

where $i \in \{0, \ldots, k\}$ is the largest integer such that $t_k^i \leq t$.

Notice that $\mathrm{(GSY)}_k$ trivially implies that

$$\frac{1}{2} \int\limits_{\Omega \backslash \Gamma} |\nabla u_k(t)|^2 \, \mathrm{d}x + \int\limits_{\Gamma} \langle g(x, \cdot), \nu_k^x(t) \rangle \, \mathrm{d}\mathcal{H}^{n-1} \leq \frac{1}{2} \int\limits_{\Omega \backslash \Gamma} |\nabla \widetilde{u}|^2 \, \mathrm{d}x + \int\limits_{\Gamma} \langle g(x, \cdot), \widetilde{\nu}^x \rangle \, \mathrm{d}\mathcal{H}^{n-1},$$

for every $\widetilde{u} \in \mathscr{A}(w_k(t))$ and for every $\widetilde{\nu} \in \mathcal{Y}(\Gamma; [0, \infty])$ with $\widetilde{\nu} \succeq \nu_k(t) \oplus \big|[\widetilde{u}] - [u_k(t)]\big|$.

*Remark* 4.3. By passing to the limit as $k \to \infty$ in $\mathrm{(IRY)}_k$, we may formally obtain the irreversibility condition for the continuous-time quasistatic evolution. (See Definition 5.1 in Section 5 below.) Unfortunately, it is not immediate to rigorously pass to the limit in $\mathrm{(IRY)}_k$: as we shall see below, in the construction of the continuous-time evolution the jumps $[u_k(t)]$ converge to $[u(t)]$ on subsequences possibly depending on $t$, thus precluding the possibility to have convergence on the same subsequence for both $[u_k(t)]$ and $[u_k(s)]$ in $\mathrm{(IRY)}_k$. For this reason, we reformulate $\mathrm{(IRY)}_k$ in a more convenient way. We start by noticing that the condition

$$V_k(t) \geq V_k(s) + \big|[u_k(t)] - [u_k(s)]\big| \quad \text{for every } s, t \in [0, T] \text{ with } s \leq t,$$

is equivalent to the system of inequalities

$$V_k(t) + [u_k(t)] \geq V_k(s) + [u_k(s)] \quad \text{for every } s, t \in [0, T] \text{ with } s \leq t, \tag{4.3}$$
$$V_k(t) - [u_k(t)] \geq V_k(s) - [u_k(s)] \quad \text{for every } s, t \in [0, T] \text{ with } s \leq t. \tag{4.4}$$

Let us notice that since $V_0 \geq \big|[u_0]\big|$ by (2.5), we have $V_k(t) + [u_k(t)] \geq 0$ and $V_k(t) - [u_k(t)] \geq 0$ for every $t \in [0, T]$. In terms of the Young measures $\nu_k$, the inequalities (4.3) and (4.4) are equivalent to stating that the functions

$$t \mapsto \nu_k(t) \oplus [u_k(t)] =: \lambda_k^{\oplus}(t) \in \mathcal{Y}(\Gamma; [0, \infty]), \tag{4.5}$$
$$t \mapsto \nu_k(t) \ominus [u_k(t)] =: \lambda_k^{\ominus}(t) \in \mathcal{Y}(\Gamma; [0, \infty]) \tag{4.6}$$

are nondecreasing with respect to $t$. Thanks to the Helly Selection Principle for Young measures (Theorem 3.17), (4.5) and (4.6) are easier to handle than $\mathrm{(IRY)}_k$, as we shall see later in Section 5.

We conclude this section with the following proposition, which shall be used to pass to the limit in $\mathrm{(GSY)}_k$ as $k \to \infty$.

**Proposition 4.4.** *Let* $w_k \rightharpoonup w$ *weakly in* $H^1(\Omega)$. *Let* $v_k \in \mathscr{A}(w_k)$ *and* $v \in H^1(\Omega \backslash \Gamma)$ *be such that* $v_k \rightharpoonup v$ *weakly in* $H^1(\Omega \backslash \Gamma)$ *and let* $\mu_k, \mu \in \mathcal{Y}(\Gamma; [0, \infty])$ *be such that* $\mu_k \to \mu$. *Let us assume that that for every* $k \in \mathbb{N}$

$$\frac{1}{2} \int\limits_{\Omega \backslash \Gamma} |\nabla v_k|^2 \, \mathrm{d}x + \int\limits_{\Gamma} \langle g(x, \cdot), \mu_k^x \rangle \, \mathrm{d}\mathcal{H}^{n-1} \leq \frac{1}{2} \int\limits_{\Omega \backslash \Gamma} |\nabla \widetilde{v}|^2 \, \mathrm{d}x + \int\limits_{\Gamma} \langle g(x, \cdot), \widetilde{\mu}_k^x \rangle \, \mathrm{d}\mathcal{H}^{n-1}, \tag{4.7}$$

*for every* $\widetilde{v} \in \mathscr{A}(w_k)$, *where* $\widetilde{\mu}_k := \mu_k \oplus \big||[\widetilde{v}] - [v_k]\big|| \in \mathcal{Y}(\Gamma; [0, \infty))$. *Then* $v \in \mathscr{A}(w)$ *and*

$$\frac{1}{2} \int\limits_{\Omega \setminus \Gamma} |\nabla v|^2 \, \mathrm{d}x + \int\limits_{\Gamma} \langle g(x, \cdot), \mu^x \rangle \, \mathrm{d}\mathcal{H}^{n-1} \leq \frac{1}{2} \int\limits_{\Omega \setminus \Gamma} |\nabla \widetilde{v}|^2 \, \mathrm{d}x + \int\limits_{\Gamma} \langle g(x, \cdot), \widetilde{\mu}^x \rangle \, \mathrm{d}\mathcal{H}^{n-1}, \qquad (4.8)$$

*for every* $\widetilde{v} \in \mathscr{A}(w)$, *where* $\widetilde{\mu} := \mu \oplus \big||[\widetilde{v}] - [v]\big|| \in \mathcal{Y}(\Gamma; [0, \infty))$.

*Proof.* By the continuity of the trace operator on $\partial_D \Omega$ with respect to the weak convergence in $H^1(\Omega \setminus \Gamma)$ we have $v \in \mathscr{A}(w)$. To prove (4.8), fix $\widetilde{v} \in \mathscr{A}(w)$. Define $\widetilde{\mu} := \mu \oplus \big||[\widetilde{v}] - [v]\big|| \in \mathcal{Y}(\Gamma; [0, \infty))$ and

$$\widetilde{v}_k := v_k + \widetilde{v} - v \in \mathscr{A}(w_k), \qquad (4.9)$$

$$\widetilde{\mu}_k := \mu_k \oplus \big||[\widetilde{v}] - [v]\big|| = \mu_k \oplus \big||[\widetilde{v}_k] - [v_k]\big||.$$

Since $v_k \rightharpoonup v$ and $\mu_k \rightharpoonup \mu$, by Remark 3.10 we have

$$\widetilde{v}_k \rightharpoonup \widetilde{v} \quad \text{weakly in } H^1(\Omega \setminus \Gamma), \qquad (4.10)$$

$$\widetilde{\mu}_k \rightharpoonup \widetilde{\mu} \quad \text{narrowly.} \qquad (4.11)$$

From (4.7) we get that

$$\frac{1}{2} \int\limits_{\Omega \setminus \Gamma} |\nabla v_k|^2 \, \mathrm{d}x + \int\limits_{\Gamma} \langle g(x, \cdot), \mu_k^x \rangle \, \mathrm{d}\mathcal{H}^{n-1} \leq \frac{1}{2} \int\limits_{\Omega \setminus \Gamma} |\nabla \widetilde{v}_k|^2 \, \mathrm{d}x + \int\limits_{\Gamma} \langle g(x, \cdot), \widetilde{\mu}_k^x \rangle \, \mathrm{d}\mathcal{H}^{n-1}. \qquad (4.12)$$

We now use a classical quadratic trick. By (4.9), we infer that

$$\begin{aligned} \frac{1}{2} \int\limits_{\Omega \setminus \Gamma} |\nabla v_k|^2 \, \mathrm{d}x - \frac{1}{2} \int\limits_{\Omega \setminus \Gamma} |\nabla \widetilde{v}_k|^2 \, \mathrm{d}x &= \frac{1}{2} \int\limits_{\Omega \setminus \Gamma} (\nabla v_k - \nabla \widetilde{v}_k) \cdot (\nabla v_k + \nabla \widetilde{v}_k) \, \mathrm{d}x \\ &= \frac{1}{2} \int\limits_{\Omega \setminus \Gamma} (\nabla v - \nabla \widetilde{v}) \cdot (2\nabla v_k + \nabla \widetilde{v} - \nabla v) \, \mathrm{d}x. \end{aligned} \qquad (4.13)$$

Thanks to (4.11) we deduce that

$$\int\limits_{\Gamma} \langle g(x, \cdot), \widetilde{\mu}_k^x \rangle \, \mathrm{d}\mathcal{H}^{n-1} \to \int\limits_{\Gamma} \langle g(x, \cdot), \widetilde{\mu}^x \rangle \, \mathrm{d}\mathcal{H}^{n-1}. \qquad (4.14)$$

Since $v_k \rightharpoonup v$ and $\mu_k \rightharpoonup \mu$, by (4.12)–(4.14) we have

$$\frac{1}{2} \int\limits_{\Omega \setminus \Gamma} (\nabla v - \nabla \widetilde{v}) \cdot (\nabla v + \nabla \widetilde{v}) \, \mathrm{d}x + \int\limits_{\Gamma} \langle g(x, \cdot), \mu^x \rangle \, \mathrm{d}\mathcal{H}^{n-1} \leq \int\limits_{\Gamma} \langle g(x, \cdot), \widetilde{\mu}^x \rangle \, \mathrm{d}\mathcal{H}^{n-1},$$

from which we easily conclude that (4.8) holds. $\qquad \square$

## 5. Quasistatic evolution in the setting of Young measures

In this section we study the continuous-time limit of the discrete evolutions $u_k(t)$ constructed in Section 4. The limit of the sequence of (Young measures concentrated on) functions $\nu_k(t)$ defined in (2.16) can only be found in the space of Young measures $\mathcal{Y}(\Gamma; [0, \infty))$. For this reason we require a definition of quasistatic evolution in a generalised sense.

**Definition 5.1.** Let $w$, $u_0$, and $V_0$ be as in (2.2)–(2.5). A *quasistatic evolution in the sense of Young measures* with initial conditions $(u_0, V_0)$ and boundary datum $w$ is a function $t \mapsto (u(t), \nu(t))$ defined in $[0, T]$ with values in $H^1(\Omega \setminus \Gamma) \times \mathcal{Y}(\Gamma; [0, \infty])$ that satisfies $u(0) = u_0$, $\nu(0) = \delta_{V_0}$, and the following conditions:

(IRY) *Irreversibility*: $\nu(t) \succeq \nu(s) \oplus \big||[u(t)] - [u(s)]\big||$ for every $s, t \in [0, T]$ with $s \leq t$.

(GSY) *Global stability*: For every $t \in [0, T]$, $u(t) \in \mathscr{A}(w(t))$ and

$$\frac{1}{2} \int\limits_{\Omega \setminus \Gamma} |\nabla u(t)|^2 \, \mathrm{d}x + \int\limits_{\Gamma} \langle g(x, \cdot), \nu^x(t) \rangle \, \mathrm{d}\mathcal{H}^{n-1} \leq \frac{1}{2} \int\limits_{\Omega \setminus \Gamma} |\nabla \widetilde{u}|^2 \, \mathrm{d}x + \int\limits_{\Gamma} \langle g(x, \cdot), \widetilde{\nu}^x \rangle \, \mathrm{d}\mathcal{H}^{n-1},$$

for every $\widetilde{u} \in \mathscr{A}(w(t))$, where $\widetilde{\nu} := \nu(t) \oplus \big||[\widetilde{u}] - [u(t)]\big|| \in \mathcal{Y}(\Gamma; [0, \infty))$.

(EBY) *Energy-dissipation balance*: For every $t \in [0, T]$

$$\frac{1}{2} \int_{\Omega \setminus \Gamma} |\nabla u(t)|^2 \, dx + \int_{\Gamma} \langle g(x, \cdot), \nu^x(t) \rangle \, d\mathcal{H}^{n-1}$$
$$= \frac{1}{2} \int_{\Omega \setminus \Gamma} |\nabla u_0|^2 \, dx + \int_{\Gamma} g(x, V_0) \, d\mathcal{H}^{n-1} + \int_0^t \langle \nabla u(s), \nabla \dot{w}(s) \rangle_{L^2} \, ds \,.$$

*Remark* 5.2. In order to recognise the connection with the classical notion of quasistatic evolution, we notice that $t \mapsto u(t)$ is a quasistatic evolution (Definition 2.6) if and only if $t \mapsto (u(t), \delta_{V_u(t)})$ is a quasistatic evolution in the sense of Young measures (Definition 5.1), where $V_u(t)$ is the function defined in (2.3). Indeed, the irreversibility condition (IRY) of Definition 5.1 automatically holds for $t \mapsto \delta_{V_u(t)}$ by definition of essential variation. Moreover, (GS) and (EB) correspond to (GSY) and (EBY), since the Young measure considered in this case is concentrated on $V_u(t)$.

*Remark* 5.3. Notice that (GSY) trivially implies that

$$\frac{1}{2} \int_{\Omega \setminus \Gamma} |\nabla u(t)|^2 \, dx + \int_{\Gamma} \langle g(x, \cdot), \nu^x(t) \rangle \, d\mathcal{H}^{n-1} \leq \frac{1}{2} \int_{\Omega \setminus \Gamma} |\nabla \widetilde{u}|^2 \, dx + \int_{\Gamma} \langle g(x, \cdot), \widetilde{\nu}^x \rangle \, d\mathcal{H}^{n-1},$$

for every $\widetilde{u} \in \mathscr{A}(w(t))$ and for every $\widetilde{\nu} \in \mathcal{Y}(\Gamma; [0, \infty])$ with $\widetilde{\nu} \succeq \nu(t) \oplus \big| [\widetilde{u}] - [u(t)] \big|$.

Moreover we underline that (IRY) is a stronger condition than the monotonicity of $t \mapsto \nu(t)$ and dictates a relationship between $\nu$ and $[u]$.

In the following theorem we prove the existence of a quasistatic evolution in the sense of of Young measures. As explained in Section 2, this result will be then improved in Section 6 by showing that the truncated Young measures $\mathcal{T}_\mu^\theta \nu(t)$ are concentrated on the function $V_u(t) \wedge \theta$ which represents the cumulation of the jumps on $\Gamma$.

**Theorem 5.4** (Existence of quasistatic evolutions in the sense of Young measures). *Assume that $g$ satisfies* (g1)–(g4) *and let $w$, $u_0$, and $V_0$ be as in* (2.2), (2.4), *and* (2.5)*. Assume that the pair* $(u_0, \delta_{V_0})$ *is globally stable, i.e.,* (2.8) *holds. Then there exists a quasistatic evolution in the sense of Young measures* $t \mapsto (u(t), \nu(t))$ *with initial conditions* $(u_0, V_0)$ *and boundary datum $w$.*

In the rest of this section, we give a proof of Theorem 5.4.

**Construction of the evolution.** Let us consider the Young measures $\nu_k(t)$ defined in (2.16). The starting point of the proof is the construction of a limit of $\nu_k(t)$ as $k \to \infty$. Since the functions $t \mapsto \nu_k(t) \in \mathcal{Y}(\Gamma; [0, \infty])$ are increasing with respect to the order $\preceq$, we can apply Theorem 3.17 to deduce that there exists a subsequence (independent of $t$ and still denoted by $\nu_k$) and an increasing function $t \mapsto \nu(t)$ from $[0, T]$ to $\mathcal{Y}(\Gamma; [0, \infty])$ such that

$$\nu_k(t) \rightharpoonup \nu(t) \quad \text{narrowly for every } t \in [0, T]. \tag{5.1}$$

Unfortunately, the convergence in (5.1) is not enough to guarantee that the irreversibility condition (IRY) holds for $\nu(t)$. In other words, it is nontrivial to pass to the limit in the discrete version of the irreversibility condition $(\text{IRY})_k$. Nonetheless, by Remark 4.3, we know that the functions $t \mapsto \lambda_k^\oplus(t)$ and $t \mapsto \lambda_k^\ominus(t)$ are increasing. Hence we can apply again Theorem 3.17 and deduce that there exists a subsequence independent of $t$ (not relabelled) and two increasing functions $t \mapsto \lambda_\oplus(t) \in \mathcal{Y}(\Gamma; [0, \infty])$ and $t \mapsto \lambda_\ominus(t) \in \mathcal{Y}(\Gamma; [0, \infty])$ such that

$$\lambda_k^\oplus(t) \rightharpoonup \lambda_\oplus(t) \quad \text{narrowly for every } t \in [0, T], \tag{5.2}$$

$$\lambda_k^\ominus(t) \rightharpoonup \lambda_\ominus(t) \quad \text{narrowly for every } t \in [0, T]. \tag{5.3}$$

The monotonicity of both the functions $\lambda_\oplus$ and $\lambda_\ominus$ encodes the irreversibility of the process in the continuous-time evolution.

We are now in a position to construct a limit of the sequence $u_k(t)$. Thanks to (4.1), we have $\|u_k(t)\|_{H^1(\Omega \setminus \Gamma)} \leq C$, where the constant $C$ is independent of $k$ and $t$. Let $t \in [0, T]$ and let $k_j(t)$ be a subsequence of $k$ such that

$$u_{k_j(t)}(t) \rightharpoonup u(t) \quad \text{weakly in } H^1(\Omega \setminus \Gamma), \tag{5.4}$$

for some function $u(t) \in H^1(\Omega \setminus \Gamma)$.

A priori, the function $u(t)$ depends on the subsequence $k_j(t)$ such that (5.4) holds. Nevertheless, we will prove below the following result.

*Remark* 5.5. Actually, we shall prove that

$$u_k(t) \to u(t) \quad \text{strongly in } H^1(\Omega \setminus \Gamma) \tag{5.5}$$

on the whole sequence (independent of $t$) found by the Helly Selection Principle (cf. (5.1)–(5.3)).

We remark that also the topology of the convergence is improved. The convergence in (5.5) will be proved later in this section by showing that the function $u(t)$ is characterised as the unique solution to a minimum problem (Proposition 5.6). The convergence with respect to the strong topology of $H^1(\Omega \setminus \Gamma)$ will be a consequence of the energy-dissipation balance (EBY).

**Proof of irreversibility.** We can now infer (IRY) from the monotonicity of the functions $\lambda_\oplus$ and $\lambda_\ominus$ obtained in (5.2) and (5.3). Indeed, from (5.4) we deduce that $[u_{k_j(t)}] \to [u(t)]$ strongly in $L^2(\Gamma)$. By (5.1) and by Remark 3.10 this implies that $\lambda^\oplus_{k_j(t)}(t) = \nu_{k_j(t)}(t) \oplus [u_{k_j(t)}(t)] \rightharpoonup \nu(t) \oplus [u(t)]$. Thus, from (5.2) we deduce that

$$\lambda_\oplus(t) = \nu(t) \oplus [u(t)], \tag{5.6}$$

and therefore that the function $t \mapsto \nu(t) \oplus [u(t)]$ is increasing. Similarly one can prove that $\lambda_\ominus(t) = \nu(t) \ominus [u(t)]$ and that $t \mapsto \nu(t) \ominus [u(t)]$ is increasing. Therefore, for every $s, t \in [0, T]$ with $s \leq t$ we have

$$\nu(t) \oplus [u(t)] \succeq \nu(s) \oplus [u(s)],$$
$$\nu(t) \ominus [u(t)] \succeq \nu(s) \ominus [u(s)].$$

It is immediate to see that the previous inequalities imply (IRY).

In order to prove (5.5), it is convenient to make the following key observations:

- the Young measures $\lambda_\oplus(t)$ and $\nu(t)$ are obtained as limits of a sequence independent of $t$;
- the jump $[u(t)]$ can be recovered just from $\lambda_\oplus(t)$ and $\nu(t)$ thanks to (5.6).

We now make precise the previous statements. We start by observing that if $x \in \Gamma$ is such that $\lambda^x_\oplus(t) = \nu^x(t) = \delta_\infty$, then $[u(t;x)]$ is not uniquely determined by (5.6). For this reason we introduce the set

$$\Gamma_N(t) := \{x \in \Gamma \ : \ \nu^x(t) \succeq \delta_{\theta(x)}\}, \tag{5.7}$$

which corresponds to the subset of $\Gamma$ where the material is completely fractured. For $\mathcal{H}^{n-1}$-a.e. $x \in \Gamma \setminus \Gamma_N(t)$ there exists a mass $m_x \in (0, 1]$ such that $F^{[-1]}_{\nu^x(t)}(m_x) \in [0, \theta(x))$, where $F^{[-1]}_{\nu^x(t)}$ is the pseudo-inverse of the cumulative distribution function $F_{\nu^x(t)}$ of $\nu^x(t)$ (cf. (3.3) and (3.4)). In particular, we have that $F^{[-1]}_{\nu^x(t)}(m_x)$ is finite. By (5.6) and by the definition of pseudo-inverse, it is easy to see that

$$F^{[-1]}_{\lambda^\oplus_\oplus(t)}(m_x) - F^{[-1]}_{\nu^x(t)}(m_x) = [u(t;x)] \quad \text{for } \mathcal{H}^{n-1}\text{-a.e. } x \in \Gamma \setminus \Gamma_N(t). \tag{5.8}$$

(We remark that, if instead $x \in \Gamma_N(t)$, it may happen that $\nu^x(t) = \delta_\infty$, and thus $F^{[-1]}_{\nu^x(t)}(m) = \infty$ for every $m \in (0, 1]$. This does not allow us to infer (5.8).) Therefore, we can define a measurable function $\gamma(t) \colon \Gamma \setminus \Gamma_N(t) \to \mathbb{R}$ by

$$\gamma(t;x) := F^{[-1]}_{\lambda^\oplus_\oplus(t)}(m_x) - F^{[-1]}_{\nu^x(t)}(m_x), \tag{5.9}$$

for $\mathcal{H}^{n-1}$-a.e. $x \in \Gamma \setminus \Gamma_N(t)$. We stress that the function $\gamma(t)$ is obtained independently of the subsequence $k_j(t)$. The proof of (5.5) will be continued after the proof of (GSY) and (EBY).

**Proof of global stability.** The global stability (GSY) directly follows from Proposition 4.4, since $u_{k_j(t)}(t)$ and $\nu_{k_j(t)}(t)$ satisfy condition $(GSY)_k$ and by (5.4) and (5.1).

In general, the function $u(t)$ is not uniquely determined by (GSY), because $u(t)$ appears both in the left-hand side and in the right-hand side of (GSY); specifically, $\widetilde{\nu}$ depends on $u(t)$. However, we have shown that the jump of $u(t)$ is given by the function $\gamma(t)$ defined in (5.9) independently of the subsequence $k_j(t)$. This allows us to prove the following result.

**Proposition 5.6.** *The function $u(t)$ obtained in (5.4) is the unique solution to the minimum problem*

$$\min_{\widetilde{u}} \left\{ \frac{1}{2} \int\limits_{\Omega \setminus \Gamma} |\nabla \widetilde{u}|^2 \, \mathrm{d}x \; : \; \widetilde{u} \in \mathscr{A}(w(t)) \text{ such that } [\widetilde{u}(x)] = \gamma(t; x) \text{ for } \mathcal{H}^{n-1}\text{-a.e. } x \in \Gamma \setminus \Gamma_N(t) \right\},$$

(5.10)

*where $\Gamma_N(t)$ is the set defined in (5.7) and $\gamma(t)$ is the function defined in (5.9).*

*Remark* 5.7. Notice that Proposition 5.6 holds true also when $\mathcal{H}^{n-1}(\Gamma \setminus \Gamma_N(t)) = 0$, i.e., when the material is completely fractured on the whole surface $\Gamma$. In this case, the competitors in (5.10) are all functions $\widetilde{u} \in \mathscr{A}(w(t))$ (without any constraint on the jump).

*Proof of Proposition 5.6.* We have already observed (see (5.8)) that $[u(t)] = \gamma(t)$ $\mathcal{H}^{n-1}$-a.e. on $\Gamma \setminus \Gamma_N(t)$. Let us fix $\widetilde{u} \in \mathscr{A}(w(t))$ such that $[\widetilde{u}] = \gamma(t) = [u(t)]$ $\mathcal{H}^{n-1}$-a.e. on $\Gamma \setminus \Gamma_N(t)$. Setting $\widetilde{\nu} := \nu(t) \oplus \big| [\widetilde{u}] - [u(t)] \big|$, by (GSY) we have

$$\frac{1}{2} \int\limits_{\Omega \setminus \Gamma} |\nabla u(t)|^2 \, \mathrm{d}x + \int\limits_{\Gamma} \langle g(x, \cdot), \nu^x(t) \rangle \, \mathrm{d}\mathcal{H}^{n-1} \leq \frac{1}{2} \int\limits_{\Omega \setminus \Gamma} |\nabla \widetilde{u}|^2 \, \mathrm{d}x + \int\limits_{\Gamma} \langle g(x, \cdot), \widetilde{\nu}^x \rangle \, \mathrm{d}\mathcal{H}^{n-1}. \tag{5.11}$$

Since $\widetilde{\nu}^x \succeq \nu^x(t) \succeq \delta_{\theta(x)}$ for $\mathcal{H}^{n-1}$-a.e. $x \in \Gamma_N(t)$ and since $g(x, \xi) = \kappa(x)$ for every $\xi \in [\theta(x), \infty]$, we deduce that

$$\int\limits_{\Gamma_N(t)} \langle g(x, \cdot), \nu^x(t) \rangle \, \mathrm{d}\mathcal{H}^{n-1} = \int\limits_{\Gamma_N(t)} \langle g(x, \cdot), \widetilde{\nu}^x \rangle \, \mathrm{d}\mathcal{H}^{n-1} = \int\limits_{\Gamma_N(t)} \kappa(x) \, \mathrm{d}\mathcal{H}^{n-1}(x).$$

Therefore (5.11) is equivalent to

$$\frac{1}{2} \int\limits_{\Omega \setminus \Gamma} |\nabla u(t)|^2 \, \mathrm{d}x + \int\limits_{\Gamma \setminus \Gamma_N(t)} \langle g(x, \cdot), \nu^x(t) \rangle \, \mathrm{d}\mathcal{H}^{n-1} \leq \frac{1}{2} \int\limits_{\Omega \setminus \Gamma} |\nabla \widetilde{u}|^2 \, \mathrm{d}x + \int\limits_{\Gamma \setminus \Gamma_N(t)} \langle g(x, \cdot), \widetilde{\nu}^x \rangle \, \mathrm{d}\mathcal{H}^{n-1}.$$

Since $[\widetilde{u}] = [u(t)]$ $\mathcal{H}^{n-1}$-a.e. on $\Gamma \setminus \Gamma_N(t)$, we have $\widetilde{\nu}^x = \nu^x(t)$ for $\mathcal{H}^{n-1}$-a.e. $x \in \Gamma \setminus \Gamma_N(t)$, hence the previous inequality reads

$$\frac{1}{2} \int\limits_{\Omega \setminus \Gamma} |\nabla u(t)|^2 \, \mathrm{d}x \leq \frac{1}{2} \int\limits_{\Omega \setminus \Gamma} |\nabla \widetilde{u}|^2 \, \mathrm{d}x.$$

This proves that $u(t)$ is a solution to (5.10).

The argument to prove uniqueness is standard: if $u_1$ and $u_2$ were two different solutions to (5.10), then $\widetilde{u} := \frac{1}{2}(u_1 + u_2)$ would be an admissible competitor; by strict convexity,

$$\frac{1}{2} \int\limits_{\Omega \setminus \Gamma} |\nabla \widetilde{u}|^2 \, \mathrm{d}x = \frac{1}{2} \int\limits_{\Omega \setminus \Gamma} \Big| \frac{\nabla u_1 + \nabla u_2}{2} \Big|^2 \, \mathrm{d}x < \frac{1}{4} \int\limits_{\Omega \setminus \Gamma} |\nabla u_1|^2 \, \mathrm{d}x + \frac{1}{4} \int\limits_{\Omega \setminus \Gamma} |\nabla u_2|^2 \, \mathrm{d}x = \frac{1}{2} \int\limits_{\Omega \setminus \Gamma} |\nabla u_1|^2 \, \mathrm{d}x.$$

This contradicts the minimality. $\qquad\square$

*Remark* 5.8. The minimum problem (5.10) is independent of the subsequence $k_j(t)$. As a consequence, we have shown that if $k_j(t)$ is such that $u_{k_j(t)} \rightharpoonup u(t)$, then $u(t)$ is the unique solution to (5.10). Thus $u(t)$ does not depend on $k_j(t)$, and this implies that

$$u_k(t) \rightharpoonup u(t) \quad \text{weakly in } H^1(\Omega \setminus \Gamma) \quad \text{for every } t \in [0, T] \tag{5.12}$$

on the whole sequence (independent of $t$) found by the Helly Selection Principle (cf. (5.1)–(5.3)). In particular, by (4.1) we have

$$\|u(t)\|_{H^1(\Omega \setminus \Gamma)} \leq C. \tag{5.13}$$

After proving the energy-dissipation balance, it will turn out that the convergence is strong in $H^1(\Omega \setminus \Gamma)$.

**Proof of energy-dissipation balance.** Before proving (EBY), we show that the function $t \mapsto u(t)$ is continuous with respect to the weak topology for almost every time. This result allows for a simple proof of the energy-dissipation balance.

**Lemma 5.9.** *There exists a countable set* $E \subset [0, T]$ *such that for every* $t \in [0, T] \setminus E$

$$u(s) \rightharpoonup u(t) \quad \text{weakly in } H^1(\Omega \setminus \Gamma), \tag{5.14}$$

$$\nu(s) \rightharpoonup \nu(t) \quad \text{narrowly in } \mathcal{Y}(\Gamma; [0, \infty]). \tag{5.15}$$

*as* $s \to t$.

*Proof.* Since the functions $t \mapsto \lambda_{\oplus}(t)$ and $t \mapsto \nu(t)$ are nondecreasing, we can find a countable set $E \subset [0, T]$ such that both $\lambda_{\oplus}$ and $\nu$ are continuous (with respect to the narrow topology) in $t$ for every $t \in [0, T] \setminus E$. (See Remark 3.16.) Thus, given $t \in [0, T] \setminus E$ and a sequence $s_k \to t$, we have

$$\lambda_{\oplus}(s_k) \rightharpoonup \lambda_{\oplus}(t), \quad \nu(s_k) \rightharpoonup \nu(t). \tag{5.16}$$

Thanks to (5.13), we can extract a subsequence (not relabelled) such that

$$u(s_k) \rightharpoonup u^* \quad \text{weakly in } H^1(\Omega \setminus \Gamma) \tag{5.17}$$

for some $u^* \in H^1(\Omega \setminus \Gamma)$. By Proposition 4.4, we infer that $u^* \in \mathscr{A}(w(t))$ and

$$\frac{1}{2} \int_{\Omega \setminus \Gamma} |\nabla u^*|^2 \, dx + \int_{\Gamma} \langle g(x, \cdot), \nu^x(t) \rangle \, d\mathcal{H}^{n-1} \leq \frac{1}{2} \int_{\Omega \setminus \Gamma} |\nabla \widetilde{u}|^2 \, dx + \int_{\Gamma} \langle g(x, \cdot), \widetilde{\nu}^x \rangle \, d\mathcal{H}^{n-1},$$

for every $\widetilde{u} \in \mathscr{A}(w(t))$, where $\widetilde{\nu} = \nu(t) \oplus \big[ |[\widetilde{u}] - [u^*]| \big]$.

On the other hand, by (5.6), we have $\lambda_{\oplus}(s_k) = \nu(s_k) \oplus [u(s_k)]$. By (5.16), (5.17), and Remark 3.10 we deduce that $\lambda_{\oplus}(t) = \nu(t) \oplus [u^*]$. Hence, by (5.9), we obtain that $[u^*](x) = \gamma(t; x)$ for $\mathcal{H}^{n-1}$-a.e. $x \in \Gamma \setminus \Gamma_N(t)$. Therefore, arguing as in the proof of Proposition 5.6, we infer that $u^*$ is a solution to the minimum problem (5.10). By uniqueness of the solution we get $u^* = u(t)$, which concludes the proof. □

*Remark* 5.10. Lemma 5.9 will be improved in Proposition 5.12 below by showing that the continuity actually holds with respect to the strong topology.

Let us now prove (EBY). We start with proving the inequality

$$\frac{1}{2} \int_{\Omega \setminus \Gamma} |\nabla u(t)|^2 \, dx + \int_{\Gamma} \langle g(x, \cdot), \nu^x(t) \rangle \, d\mathcal{H}^{n-1}$$
$$\leq \frac{1}{2} \int_{\Omega \setminus \Gamma} |\nabla u_0|^2 \, dx + \int_{\Gamma} g(x, V_0) \, d\mathcal{H}^{n-1} + \int_0^t \langle \nabla u(s), \nabla \dot{w}(s) \rangle_{L^2} \, ds. \tag{5.18}$$

By (5.1), (5.12), and by (EIY)$_k$, for every $t \in [0, T]$ we have

$$\frac{1}{2} \int_{\Omega \setminus \Gamma} |\nabla u(t)|^2 \, dx + \int_{\Gamma} \langle g(x, \cdot), \nu^x(t) \rangle \, d\mathcal{H}^{n-1}$$
$$\leq \liminf_{k \to \infty} \left[ \frac{1}{2} \int_{\Omega \setminus \Gamma} |\nabla u_k(t)|^2 \, dx + \int_{\Gamma} \langle g(x, \cdot), \nu_k^x(t) \rangle \, d\mathcal{H}^{n-1} \right] \tag{5.19}$$
$$\leq \frac{1}{2} \int_{\Omega \setminus \Gamma} |\nabla u_0|^2 \, dx + \int_{\Gamma} g(x, V_0) \, d\mathcal{H}^{n-1} + \limsup_{k \to \infty} \int_0^{t_k^i} \langle \nabla u_k(s), \nabla \dot{w}(s) \rangle_{L^2} \, ds,$$

where $i \in \{0, \ldots, k\}$ is the largest integer such that $t_k^i \leq t$. Thanks to (5.12) we know that

$$\langle \nabla u_k(s), \nabla \dot{w}(s) \rangle_{L^2} \to \langle \nabla u(s), \nabla \dot{w}(s) \rangle_{L^2} \quad \text{for every } s \in [0, t].$$

Moreover, from (4.1) we deduce that

$$\langle \nabla u_k(s), \nabla \dot{w}(s) \rangle_{L^2} \leq \|\nabla u_k(s)\|_{L^2} \|\nabla \dot{w}(s)\|_{L^2} \leq C \|\nabla \dot{w}(s)\|_{L^2},$$

for every $s \in [0, T]$. By our assumption (2.2) on $w$, the function $t \mapsto \nabla \dot{w}(t)$ is $L^1([0, T]; L^2(\Omega \setminus \Gamma))$, so we can apply the Dominated Convergence Theorem to infer that

$$\limsup_{k \to \infty} \int_0^{t_k^i} \langle \nabla u_k(s), \nabla \dot{w}(s) \rangle_{L^2} \, \mathrm{d}s = \lim_{k \to \infty} \int_0^t \langle \nabla u_k(s), \nabla \dot{w}(s) \rangle_{L^2} \, \mathrm{d}s = \int_0^t \langle \nabla u(s), \nabla \dot{w}(s) \rangle_{L^2} \, \mathrm{d}s. \quad (5.20)$$

Together with (5.19), the previous inequality yields (5.18).

We now exploit the global stability to prove, for a fixed $t \in [0, T]$, the opposite inequality

$$\begin{aligned}
&\frac{1}{2} \int_{\Omega \setminus \Gamma} |\nabla u(t)|^2 \, \mathrm{d}x + \int_\Gamma \langle g(x, \cdot), \nu^x(t) \rangle \, \mathrm{d}\mathcal{H}^{n-1} \\
&\geq \frac{1}{2} \int_{\Omega \setminus \Gamma} |\nabla u_0|^2 \, \mathrm{d}x + \int_\Gamma g(x, V_0) \, \mathrm{d}\mathcal{H}^{n-1} + \int_0^t \langle \nabla u(s), \nabla \dot{w}(s) \rangle_{L^2} \, \mathrm{d}s.
\end{aligned} \quad (5.21)$$

For every $k \in \mathbb{N}$, let us consider the subdivision of the time interval $[0, t]$ given by the $k+1$ equispaced nodes

$$s_k^h := \tfrac{h}{k} t \quad \text{for } h = 0, \dots, k.$$

Let $h \in \{1, \dots, k\}$. By the irreversibility condition (IRY), we have $\nu(s_k^h) \succeq \nu(s_k^{h-1}) \oplus [[u(s_k^h)] - [u(s_k^{h-1})]] =: \tilde{\nu}_h$. Since $u(s_k^h) - w(s_k^h) + w(s_k^{h-1}) \in \mathscr{A}(w(s_k^{h-1}))$, by (GSY) we obtain

$$\begin{aligned}
&\frac{1}{2} \int_{\Omega \setminus \Gamma} |\nabla u(s_k^{h-1})|^2 \, \mathrm{d}x + \int_\Gamma \langle g(x, \cdot), \nu^x(s_k^{h-1}) \rangle \, \mathrm{d}\mathcal{H}^{n-1} \\
&\leq \frac{1}{2} \int_{\Omega \setminus \Gamma} |\nabla u(s_k^h)|^2 \, \mathrm{d}x + \int_\Gamma \langle g(x, \cdot), \tilde{\nu}_h^x \rangle \, \mathrm{d}\mathcal{H}^{n-1} \\
&\quad - \int_{\Omega \setminus \Gamma} \nabla u(s_k^h) \cdot (\nabla w(s_k^h) - \nabla w(s_k^{h-1})) \, \mathrm{d}x + \frac{1}{2} \int_{\Omega \setminus \Gamma} |\nabla w(s_k^h) - \nabla w(s_k^{h-1})|^2 \, \mathrm{d}x \\
&\leq \frac{1}{2} \int_{\Omega \setminus \Gamma} |\nabla u(s_k^h)|^2 \, \mathrm{d}x + \int_\Gamma \langle g(x, \cdot), \nu^x(s_k^h) \rangle \, \mathrm{d}\mathcal{H}^{n-1} \\
&\quad - \int_{s_k^{h-1}}^{s_k^h} \langle \nabla \bar{u}^k(s), \nabla \dot{w}(s) \rangle_{L^2} \, \mathrm{d}s + \frac{1}{2} \left( \int_{s_k^{h-1}}^{s_k^h} \|\nabla \dot{w}(s)\|_{L^2} \, \mathrm{d}s \right)^2,
\end{aligned} \quad (5.22)$$

where

$$\bar{u}^k(s) := u(s_k^h) \quad \text{for every } s \in (s_k^{h-1}, s_k^h].$$

Summing up the inequalities given by (5.22) for $h = 1, \dots, k$, we get

$$\begin{aligned}
&\frac{1}{2} \int_{\Omega \setminus \Gamma} |\nabla u(t)|^2 \, \mathrm{d}x + \int_\Gamma \langle g(x, \cdot), \nu^x(t) \rangle \, \mathrm{d}\mathcal{H}^{n-1} \\
&\geq \frac{1}{2} \int_{\Omega \setminus \Gamma} |\nabla u_0|^2 \, \mathrm{d}x + \int_\Gamma g(x, V_0) \, \mathrm{d}\mathcal{H}^{n-1} + \int_0^t \langle \nabla \bar{u}_k(s), \nabla \dot{w}(s) \rangle_{L^2} \, \mathrm{d}s - \bar{\eta}_k,
\end{aligned}$$

where

$$\bar{\eta}_k := \frac{1}{2} \left( \max_{1 \leq h \leq k} \int_{s_k^{h-1}}^{s_k^h} \|\nabla \dot{w}(s)\|_{L^2} \, \mathrm{d}s \right) \left( \int_0^T \|\nabla \dot{w}(s)\|_{L^2} \, \mathrm{d}s \right).$$

In order to infer (5.21), we notice that by Lemma 5.9 we have $\bar{u}^k(s) \rightharpoonup u(s)$ for almost every $s \in [0, t]$, and therefore

$$\lim_{k \to \infty} \int_0^t \langle \nabla \bar{u}^k(s), \nabla \dot{w}(s) \rangle_{L^2} \, \mathrm{d}s = \int_0^t \langle \nabla u(s), \nabla \dot{w}(s) \rangle_{L^2} \, \mathrm{d}s,$$

by the Dominated Convergence Theorem. This concludes the proof of (EBY) and of Theorem 5.4.

**Approximation of the evolution and continuity for almost every time.** Thanks to (EBY), we prove the convergence of the approximating evolutions (5.5) and we improve Lemma 5.9.

**Proposition 5.11.** *We have*

$$u_k(t) \to u(t) \quad \text{strongly in } H^1(\Omega \setminus \Gamma)$$

*on the whole sequence (independent of $t$) such that (5.1)–(5.3) hold.*

*Proof.* By (5.1) and (5.12), for every $t \in [0, T]$ we have

$$\frac{1}{2} \int_{\Omega \setminus \Gamma} |\nabla u(t)|^2 \, dx + \int_\Gamma \langle g(x, \cdot), \nu^x(t) \rangle \, d\mathcal{H}^{n-1} \le \liminf_{k \to \infty} \left[ \frac{1}{2} \int_\Omega |\nabla u_k(t)|^2 \, dx + \int_\Gamma \langle g(x, \cdot), \nu_k^x(t) \rangle \, d\mathcal{H}^{n-1} \right]. \tag{5.23}$$

On the other hand, by (5.20), (EBY), and $(\text{EIY})_k$ we get

$$\begin{aligned}
\limsup_{k \to \infty} & \left[ \frac{1}{2} \int_\Omega |\nabla u_k(t)|^2 \, dx + \int_\Gamma \langle g(x, \cdot), \nu_k^x(t) \rangle \, d\mathcal{H}^{n-1} \right] \\
& \le \frac{1}{2} \int_{\Omega \setminus \Gamma} |\nabla u_0|^2 \, dx + \int_\Gamma g(x, V_0) \, d\mathcal{H}^{n-1} + \int_0^t \langle \nabla u(s), \nabla \dot{w}(s) \rangle_{L^2} \, ds \\
& = \frac{1}{2} \int_{\Omega \setminus \Gamma} |\nabla u(t)|^2 \, dx + \int_\Gamma \langle g(x, \cdot), \nu^x(t) \rangle \, d\mathcal{H}^{n-1}.
\end{aligned} \tag{5.24}$$

Thus all inequalities in (5.23) and (5.24) are equalities. Since

$$\int_\Gamma \langle g(x, \cdot), \nu_k^x(t) \rangle \, d\mathcal{H}^{n-1} \to \int_\Gamma \langle g(x, \cdot), \nu^x(t) \rangle \, d\mathcal{H}^{n-1},$$

we have $\|\nabla u_k(t)\|_{L^2} \to \|\nabla u(t)\|_{L^2}$. Thanks to (5.12), this concludes the proof. $\square$

**Proposition 5.12.** *There exists a countable set $E \subset [0, T]$ such that for every $t \in [0, T] \setminus E$*

$$u(s) \to u(t) \quad \text{strongly in } H^1(\Omega \setminus \Gamma), \tag{5.25}$$

$$\nu(s) \rightharpoonup \nu(t) \quad \text{narrowly in } \mathcal{Y}(\Gamma; [0, \infty]). \tag{5.26}$$

*as $s \to t$.*

*Proof.* By (EBY) we have for every $s, t \in [0, T]$

$$\begin{aligned}
\frac{1}{2} & \int_{\Omega \setminus \Gamma} |\nabla u(t)|^2 \, dx + \int_\Gamma \langle g(x, \cdot), \nu^x(t) \rangle \, d\mathcal{H}^{n-1} \\
& = \frac{1}{2} \int_{\Omega \setminus \Gamma} |\nabla u(s)|^2 \, dx + \int_\Gamma \langle g(x, \cdot), \nu^x(s) \rangle \, d\mathcal{H}^{n-1} + \int_s^t \langle \nabla u(r), \nabla \dot{w}(r) \rangle_{L^2} \, dr.
\end{aligned}$$

Thus, if $t$ is a continuity point for the nondecreasing function $s \mapsto \nu(s)$, we have $\|\nabla u(s)\|_{L^2} \to \|\nabla u(t)\|_{L^2}$ as $s \to t$, since $r \mapsto \langle \nabla u(r), \nabla w(r) \rangle_{L^2}$ is in $L^1([0, T])$ by (2.2) and (5.13). By Lemma 5.9, this gives the desired convergence. $\square$

## 6. Proof of the main result

This section is devoted to the proof of Theorem 2.9. Besides, we also give a proof of Proposition 2.7 and of Proposition 2.8 regarding the strong formulation of the quasistatic evolution.

In Section 5 we have shown the existence of a quasistatic evolution $(u(t), \nu(t))$ in the sense of Young measures. We will now exploit the concavity of $g(x, \cdot)$ to prove that the very same displacement $t \mapsto u(t)$ found in Section 5 is also a quasistatic evolution in the sense of Definition 2.6. We recall that $g(x, \cdot)$ is strictly increasing in the interval $[0, \theta(x)]$, $\theta(x)$ is the threshold defined in (2.6). This allows us to prove that the Young measure $\nu(t)$ truncated by $\theta$ (see (3.10) for the definition) is actually concentrated on $V_u(t) \wedge \theta$, i.e., $V_u(t) \wedge \theta$ is the limit of $V_k(t) \wedge \theta$.

*Proof of Theorem 2.9.* By Theorem 5.4 and Proposition 5.11, we know that there exists a quasistatic evolution in the sense of Young measures $t \mapsto (u(t), \nu(t))$ such that, for every $t \in [0, T]$, we have (2.12) and

$$\delta_{V_k(t)} = \nu_k(t) \rightharpoonup \nu(t) \quad \text{in } \mathcal{Y}(\Gamma; [0, \infty]), \tag{6.1}$$

up to a subsequence independent of $t$ (not relabelled).

In order to prove (GS), we first prove that

$$\nu(t) \succeq \delta_{V_u(t)} \quad \text{for every } t \in [0, T]. \tag{6.2}$$

By definition of $V_u(t)$ and Remark 3.15, it is enough to show that for any partition $P$ of $[0, t]$, $P = \{0 = s_0 < s_1 < \cdots < s_{j-1} < s_j = t\}$, we have

$$\nu(t) \succeq \delta_{V^P(t)}, \tag{6.3}$$

where

$$V^P(t) := V_0 + \sum_{i=1}^{j} \big| [u(s_i)] - [u(s_{i-1})] \big|.$$

The irreversibility condition (IRY) satisfied by $s \mapsto \nu(s)$ yields

$$\nu(s_i) \succeq \nu(s_{i-1}) \oplus \big| [u(s_i)] - [u(s_{i-1})] \big| \quad \text{for } i = 1, \ldots, j. \tag{6.4}$$

Employing (6.4) inductively, we obtain the chain of inequalities

$$\nu(t) = \nu(s_j) \succeq \nu(s_{j-1}) \oplus \big| [u(s_j)] - [u(s_{j-1})] \big|$$

$$\succeq \nu(s_{j-2}) \oplus \big( \big| [u(s_{j-1})] - [u(s_{j-2})] \big| + \big| [u(s_j)] - [u(s_{j-1})] \big| \big) \succeq \ldots$$

$$\succeq \nu(s_1) \oplus \sum_{i=2}^{j} \big| [u(s_i)] - [u(s_{i-1})] \big|$$

$$\succeq \nu(0) \oplus \sum_{i=1}^{j} \big| [u(s_i)] - [u(s_{i-1})] \big| = \delta_{V^P(t)},$$

and thus (6.2) holds true.

Recalling the definition of cumulative distribution function (3.3), we have $F_{\delta_{V_u(t;x)}}(\xi) = 0$ for $\xi < V_u(t; x)$. Thus, by (ii) in Definition 3.12, we deduce that

$$\operatorname{supp} \nu^x(t) \subset [V_u(t; x), \infty] \tag{6.5}$$

for every $t \in [0, T]$ and for $\mathcal{H}^{n-1}$-a.e. $x \in \Gamma$.

We are now in a position to prove that $t \mapsto u(t)$ satisfies the global stability condition (GS). We start by fixing $t \in [0, T]$ and $\widetilde{u} \in \mathscr{A}(w(t))$, and by setting

$$\widetilde{\nu} := \nu(t) \oplus \big| [\widetilde{u}] - [u(t)] \big|. \tag{6.6}$$

Condition (GSY) for $t \mapsto (u(t), \nu(t))$ gives

$$\frac{1}{2} \int_{\Omega \setminus \Gamma} |\nabla u(t)|^2 \, dx + \int_{\Gamma} \langle g(x, \cdot), \nu^x(t) \rangle \, d\mathcal{H}^{n-1} \leq \frac{1}{2} \int_{\Omega \setminus \Gamma} |\nabla \widetilde{u}|^2 \, dx + \int_{\Gamma} \langle g(x, \cdot), \widetilde{\nu}^x \rangle \, d\mathcal{H}^{n-1},$$

and thus (GS) follows if we show that

$$\int_{\Gamma} \Big( \langle g(x, \cdot), \widetilde{\nu}^x \rangle - \langle g(x, \cdot), \nu^x(t) \rangle \Big) \, d\mathcal{H}^{n-1} \leq \int_{\Gamma} \Big( g\big(x, V_u(t) + |[\widetilde{u}] - [u(t)]| \big) - g\big(x, V_u(t)\big) \Big) \, d\mathcal{H}^{n-1}. \tag{6.7}$$

In order to prove (6.7), notice that by (6.5) and (6.6) we have

$$\langle g(x, \cdot), \widetilde{\nu}^x \rangle - \langle g(x, \cdot), \nu^x(t) \rangle = \int_{[0, \infty]} \Big( g\big(x, \xi + |[\widetilde{u}(x)] - [u(t;x)]| \big) - g(x, \xi) \Big) \nu^x(t)(d\xi)$$

$$= \int_{[V_u(t;x), \infty]} \Big( g\big(x, \xi + |[\widetilde{u}(x)] - [u(t;x)]| \big) - g(x, \xi) \Big) \nu^x(t)(d\xi), \tag{6.8}$$

for $\mathcal{H}^{n-1}$-a.e. $x \in \Gamma$. Since $g(x, \cdot)$ is a concave function, for every $\xi \geq V_u(t; x)$ it holds

$$g\big(x, \xi + |[\widetilde{u}(x)] - [u(t;x)]| \big) - g(x, \xi) \leq g\big(x, V_u(t;x) + |[\widetilde{u}(x)] - [u(t;x)]| \big) - g\big(x, V_u(t;x)\big). \tag{6.9}$$

Let us observe that the right hand side in the inequality above does not depend on $\xi$. Therefore, by (6.8), (6.9), and recalling that $\nu^x(t)$ is a probability measure for $\mathcal{H}^{n-1}$-a.e. $x \in \Gamma$, we deduce (6.7). This completes the proof of (GS).

Let us now prove that $t \mapsto u(t)$ satisfies (EB). Arguing as in the proof of (5.21), using (GS) it is possible to see that

$$
\begin{aligned}
\frac{1}{2} \int_{\Omega \setminus \Gamma} |\nabla u(t)|^2 \, \mathrm{d}x &+ \int_{\Gamma} g\big(x, V_u(t)\big) \, \mathrm{d}\mathcal{H}^{n-1} \\
&\geq \frac{1}{2} \int_{\Omega \setminus \Gamma} |\nabla u_0|^2 \, \mathrm{d}x + \int_{\Gamma} g\big(x, V_0\big) \, \mathrm{d}\mathcal{H}^{n-1} + \int_0^t \langle \nabla u(s), \nabla \dot{w}(s) \rangle_{L^2} \, \mathrm{d}s \, .
\end{aligned}
$$

On the other hand, the opposite inequality follows immediately from (EBY) since by (6.2) we have

$$
\int_{\Gamma} g\big(x, V_u(t)\big) \, \mathrm{d}\mathcal{H}^{n-1} \leq \int_{\Gamma} \langle g(x, \cdot), \nu^x(t) \rangle \, \mathrm{d}\mathcal{H}^{n-1} \, .
$$

Therefore, $t \mapsto u(t)$ is a quasistatic evolution in the sense of Definition 2.6.

We now claim that the truncation $\mathcal{T}_\#^\theta \nu(t)$ (see (3.10) for the definition) is concentrated on $V_u(t) \wedge \theta$. To this end, we compare (EB) and (EBY), and deduce that for every $t \in [0, T]$

$$
\int_{\Gamma} g\big(x, V_u(t)\big) \, \mathrm{d}\mathcal{H}^{n-1} = \int_{\Gamma} \langle g(x, \cdot), \nu^x(t) \rangle \, \mathrm{d}\mathcal{H}^{n-1} \, . \tag{6.10}
$$

Since by (6.2) and Definition 3.12 we have $g\big(x, V_u(t; x)\big) \leq \langle g(x, \cdot), \nu^x(t) \rangle$, equality (6.10) implies that

$$
g\big(x, V_u(t; x)\big) = \langle g(x, \cdot), \nu^x(t) \rangle \tag{6.11}
$$

for $\mathcal{H}^{n-1}$-a.e. $x \in \Gamma$. Let us now fix $t$ and let $x$ be such that (6.5) holds. To prove the claim, we need to show that if $V_u(t; x) < \theta(x)$, then $\nu^x(t)\big((V_u(t; x), \infty]\big) = 0$. Let us assume, on the contrary, that $\nu^x(t)\big((V_u(t; x), \infty]\big) = c \in (0, 1]$. By (6.5) we know that

$$
\langle g(x, \cdot), \nu^x(t) \rangle = g(x, V_u(t; x))(1 - c) + \int_{(V_u(t; x), \infty]} g(x, \xi) \, \nu^x(t)(\mathrm{d}\xi) \, ,
$$

and thus

$$
\langle g(x, \cdot), \nu^x(t) \rangle - g\big(x, V_u(t; x)\big) = \int_{(V_u(t; x), \infty]} \big( g(x, \xi) - g(x, V_u(t; x)) \big) \, \nu^x(t)(\mathrm{d}\xi) \, . \tag{6.12}
$$

Since $g(x, \cdot)$ is strictly increasing in $[0, \theta(x)]$ and $\nu^x(t)\big((V_u(t; x), \infty]\big) > 0$, we get that the right-hand side in (6.12) is strictly positive. This contradicts (6.11), and therefore we have proved that $\mathcal{T}_\#^\theta \nu(t)$ is concentrated on $V_u(t) \wedge \theta$.

Eventually, using also (6.1) and Remark 3.11, we deduce that

$$
\delta_{V_k(t) \wedge \theta} = \mathcal{T}_\#^\theta \nu_k(t) \rightharpoonup \mathcal{T}_\#^\theta \nu(t) = \delta_{V_u(t) \wedge \theta} \quad \text{in } \mathcal{Y}(\Gamma; [0, \infty]) \, . \tag{6.13}
$$

By Proposition 3.6, (6.13) is equivalent to (2.13).

As for the proof of (2.14) and (2.15), we notice that by Proposition 5.12 there exists a set $E$, at most countable, such that we have (2.14) and $\nu(s) \rightharpoonup \nu(t)$ in $\mathcal{Y}(\Gamma; [0, \infty])$, for $t \in [0, T] \setminus E$ and $s \to t$. The convergence in (2.15) then follows with an argument analogous to the one used to show (2.13).

This concludes the proof. $\qquad \square$

*Remark* 6.1. In the proof of Theorem 2.9, we have shown that $\mathcal{T}_\#^\theta \nu(t) = \delta_{V_u(t) \wedge \theta}$. In particular, this allows us to rewrite the set $\Gamma_N(t)$ introduced in (5.7) (corresponding to the part of $\Gamma$ where the material is completely fractured) in terms of the variation of the jumps $V_u(t)$ and the threshold $\theta$. Namely, we have

$$
\Gamma_N(t) = \{ x \in \Gamma \ : \ V_u(t; x) \geq \theta(x) \} \, .
$$

We now give the proof of the results concerning the strong formulation of the quasistatic evolution discussed in Section 2. The derivation of the Euler-Lagrange conditions follows by standard arguments illustrated below.

*Proof of Proposition 2.7.* Let consider the set $\Gamma_N(t) = \{x \in \Gamma \; : \; V_u(t;x) \geq \theta(x)\}$. Let $\psi \in H^1(\Omega \setminus \Gamma)$ with $\psi = 0$ on $\partial_D \Omega$ and let $\varepsilon \in \mathbb{R}$. Since

$$\int\limits_{\Gamma_N(t)} g\big(x, V_u(t)\big) \, \mathrm{d}\mathcal{H}^{n-1} = \int\limits_{\Gamma_N(t)} \kappa(x) \, \mathrm{d}\mathcal{H}^{n-1} = \int\limits_{\Gamma_N(t)} g\big(x, V_u(t) + |\varepsilon[\psi]|\big) \, \mathrm{d}\mathcal{H}^{n-1}$$

and $u(t) + \varepsilon \psi \in \mathscr{A}(w(t))$, by (GS) we have

$$\frac{1}{2} \int\limits_{\Omega \setminus \Gamma} |\nabla u(t)|^2 \, \mathrm{d}x + \int\limits_{\Gamma \setminus \Gamma_N(t)} g\big(x, V_u(t)\big) \, \mathrm{d}\mathcal{H}^{n-1}$$
$$\leq \frac{1}{2} \int\limits_{\Omega \setminus \Gamma} |\nabla u(t) + \varepsilon \nabla \psi|^2 \, \mathrm{d}x + \int\limits_{\Gamma \setminus \Gamma_N(t)} g\big(x, V_u(t) + |\varepsilon[\psi]|\big) \, \mathrm{d}\mathcal{H}^{n-1}.$$

Since $g$ is of class $\mathcal{C}^1$, deriving the previous inequality with respect to $\varepsilon$ for $\varepsilon > 0$ and $\varepsilon < 0$, we get

$$-\int\limits_{\Gamma \setminus \Gamma_N(t)} g'\big(x, V_u(t)\big) \big|[\psi]\big| \, \mathrm{d}\mathcal{H}^{n-1} \leq \int\limits_{\Omega \setminus \Gamma} \nabla u(t) \cdot \nabla \psi \, \mathrm{d}x \leq \int\limits_{\Gamma \setminus \Gamma_N(t)} g'\big(x, V_u(t)\big) \big|[\psi]\big| \, \mathrm{d}\mathcal{H}^{n-1}.$$

Using the fact that $g'(x, \xi) = 0$ for $\xi \geq \theta(x)$, we also get

$$-\int\limits_{\Gamma} g'\big(x, V_u(t)\big) \big|[\psi]\big| \, \mathrm{d}\mathcal{H}^{n-1} \leq \int\limits_{\Omega \setminus \Gamma} \nabla u(t) \cdot \nabla \psi \, \mathrm{d}x \leq \int\limits_{\Gamma} g'\big(x, V_u(t)\big) \big|[\psi]\big| \, \mathrm{d}\mathcal{H}^{n-1}. \qquad (6.14)$$

By (6.14) for arbitrary $\psi \in H^1(\Omega)$ with $\psi = 0$ on $\partial_D \Omega$ and $\psi = 0$ in $\Omega^-$, we infer that $\Delta u(t) = 0$ in $\Omega^+$ and $\partial_\nu u(t) = 0$ in $H^{-\frac{1}{2}}(\partial_N \Omega \cap \partial \Omega^+)$. With similar arguments, we obtain analogous properties in $\Omega^-$ and we eventually deduce (i).

Let us prove (ii). Since $\nu_\Gamma$ is chosen in such a way that it coincides with the outer normal to $\partial \Omega^-$, by definition of normal derivative of the function $u(t)^+ = u(t)|_{\Omega^+}$ on $\Gamma$ we have that $\partial_\nu u(t)^+ \in H^{-\frac{1}{2}}(\Gamma)$ is given by

$$\langle \partial_\nu u(t)^+, \psi^+ \rangle = -\int\limits_{\Omega^+} \nabla u(t) \cdot \nabla \psi^+ \, \mathrm{d}x,$$

for every $\psi^+ \in H^1(\Omega^+)$ with $\psi^+ = 0$ on $\partial_D \Omega \cap \partial \Omega^+$. Similarly, the normal derivative $\partial_\nu u(t)^- \in H^{-\frac{1}{2}}(\Gamma)$ is given by

$$\langle \partial_\nu u(t)^-, \psi^- \rangle = \int\limits_{\Omega^-} \nabla u(t) \cdot \nabla \psi^- \, \mathrm{d}x,$$

for every $\psi^- \in H^1(\Omega^-)$ with $\psi^- = 0$ on $\partial_D \Omega \cap \partial \Omega^-$. Hence, by testing (6.14) with functions $\psi \in H^1(\Omega \setminus \Gamma)$ with $\psi = 0$ on $\partial_D \Omega$ and $[\psi] = 0$ on $\Gamma$, we infer

$$-\langle \partial_\nu u(t)^+, \psi \rangle + \langle \partial_\nu u(t)^-, \psi \rangle = 0,$$

which implies (ii) by the arbitrariness of $\psi$.

In order to prove (iii), we note that since $g'(x, \xi) \leq g'(x, 0)$ for $\mathcal{H}^{n-1}$-a.e. $x \in \Gamma$ and for every $\xi \in [0, \infty]$, by inequality (6.14) we get

$$\big| \langle \partial_\nu u(t), [\psi] \rangle \big| \leq \|g'(\cdot, 0)\|_{L^\infty} \|[\psi]\|_{L^1},$$

for every $\psi \in H^1(\Omega \setminus \Gamma)$ with $\psi = 0$ on $\partial_D \Omega$. Thus $\partial_\nu u(t)$ is a linear and continuous operator on the space $X := \{[\psi] \; : \; \psi \in H^1(\Omega \setminus \Gamma)$ such that $\psi = 0$ on $\partial_D \Omega\}$. By density of $X$ in $L^1(\Gamma)$, this implies that $\partial_\nu u(t)$ can be extended to a linear and continuous operator on $L^1(\Gamma)$, and hence $\partial_\nu u(t) \in L^\infty(\Gamma)$. From (6.14) we deduce that

$$-\int\limits_{\Gamma} g'\big(x, V_u(t)\big) |z| \, \mathrm{d}\mathcal{H}^{n-1} \leq -\int\limits_{\Gamma} \partial_\nu u(t) z \, \mathrm{d}\mathcal{H}^{n-1} \leq \int\limits_{\Gamma} g'\big(x, V_u(t)\big) |z| \, \mathrm{d}\mathcal{H}^{n-1},$$

for every $z \in L^1(\Gamma)$. This concludes the proof of (iii).                    $\square$

In order to give a proof of Proposition 2.8, we need to prove the following lemma regarding the differentiability in time of the essential variation of a function that is absolutely continuous in time with values in $L^2(\Gamma)$.

**Lemma 6.2.** *Let* $\gamma \in AC([0, T]; L^2(\Gamma))$. *Then* $\mathrm{ess\,Var}(\gamma; 0, \cdot) \in AC([0, T]; L^2(\Gamma))$ *and*

$$\lim_{s \to t} \frac{\mathrm{ess\,Var}(\gamma; s, t)}{t - s}(x) = |\dot{\gamma}(t; x)| \quad \text{for } \mathcal{H}^{n-1}\text{-a.e. } x \in \Gamma \text{ and for a.e. } t \in [0, T], \tag{6.15}$$

*where the limit and the derivative* $\dot{\gamma}$ *are defined with respect to the strong topology in* $L^2(\Gamma)$.

*Proof.* We fix $s, t \in [0, T]$ with $s < t$ and we consider a partition of the interval $[s, t]$, namely $s = s_0 < \cdots < s_j = t$. By the absolute continuity of $\gamma$, for every $i = 1, \ldots, j$ we have

$$|\gamma(s_i; x) - \gamma(s_{i-1}; x)| = \left| \int_{s_{i-1}}^{s_i} \dot{\gamma}(\tau; x) \, \mathrm{d}\tau \right| \leq \int_{s_{i-1}}^{s_i} |\dot{\gamma}(\tau; x)| \, \mathrm{d}\tau \quad \text{for } \mathcal{H}^{n-1}\text{-a.e. } x \in \Gamma,$$

where the integrals are Bochner integrals and $\dot{\gamma}(\tau)$ is the derivative in $L^2(\Gamma)$ of $\gamma(\tau)$. Summing up the previous inequalities for $i = 1, \ldots, j$, we obtain

$$\sum_{i=1}^{j} |\gamma(s_i; x) - \gamma(s_{i-1}; x)| \leq \int_{s}^{t} |\dot{\gamma}(\tau; x)| \, \mathrm{d}\tau \quad \text{for } \mathcal{H}^{n-1}\text{-a.e. } x \in \Gamma. \tag{6.16}$$

By Definition 2.3, (6.16) implies that

$$\mathrm{ess\,Var}(\gamma; s, t)(x) \leq \int_{s}^{t} |\dot{\gamma}(\tau; x)| \, \mathrm{d}\tau \quad \text{for } \mathcal{H}^{n-1}\text{-a.e. } x \in \Gamma. \tag{6.17}$$

In particular, choosing $s = 0$ in (6.17) we deduce that $\mathrm{ess\,Var}(\gamma; 0, t)$ belongs to $L^2(\Gamma)$, for every $t \in [0, T]$. By taking the $L^2$ norm in (6.17) we infer

$$\|\mathrm{ess\,Var}(\gamma; s, t)\|_{L^2} \leq \int_{s}^{t} \|\dot{\gamma}(\tau)\|_{L^2} \, \mathrm{d}\tau.$$

Since the function $\tau \mapsto \|\dot{\gamma}(\tau)\|_{L^2}$ belongs to $L^1([0, T]; \mathbb{R})$, we conclude that $\mathrm{ess\,Var}(\gamma; 0, \cdot) \in AC([0, T]; L^2(\Gamma))$.

We now compute the derivative of $\mathrm{ess\,Var}(\gamma; 0, \cdot)$. Since $\frac{1}{t-s} \int_{s}^{t} |\dot{\gamma}(\tau)| \, \mathrm{d}\tau \to |\dot{\gamma}(t)|$ strongly in $L^2(\Gamma)$ as $s \to t$, dividing all terms in (6.17) by $t - s$ and letting $s \to t$ we deduce that

$$\lim_{s \to t} \frac{\mathrm{ess\,Var}(\gamma; s, t)}{t - s}(x) \leq |\dot{\gamma}(t; x)|$$

for $\mathcal{H}^{n-1}$-a.e. $x \in \Gamma$. On the other hand, since $\{s, t\}$ is a particular partition of the interval $[s, t]$, by definition of essential variation we have

$$|\gamma(t; x) - \gamma(s; x)| \leq \mathrm{ess\,Var}(\gamma; s, t)(x),$$

for $\mathcal{H}^{n-1}$-a.e. $x \in \Gamma$. Dividing by $t - s$ and letting $s \to t$ in the inequality above, we obtain (6.15). $\square$

We are now in a position to prove Proposition 2.8.

*Proof of Proposition 2.8.* Since by assumption $u \in AC([0, T]; H^1(\Omega \setminus \Gamma))$, we have

$$\frac{\mathrm{d}}{\mathrm{d}t} \int_{\Omega \setminus \Gamma} |\nabla u(t)|^2 \, \mathrm{d}x = \int_{\Omega \setminus \Gamma} \nabla u(t) \cdot \nabla \dot{u}(t) \, \mathrm{d}x. \tag{6.18}$$

Moreover we claim that

$$\frac{\mathrm{d}}{\mathrm{d}t} \int_{\Gamma} g(x, V_u(t)) \, \mathrm{d}\mathcal{H}^{n-1} = \int_{\Gamma} g'(x, V_u(t)) |\dot{u}(t)| \, \mathrm{d}\mathcal{H}^{n-1}. \tag{6.19}$$

Let us prove (6.19). The absolute continuity of $u$ implies that $[u] \in AC([0, T]; L^2(\Gamma))$. Let us consider the set $\Gamma_N(0) = \{x \in \Gamma : V_0(x) \geq \theta(x)\}$. Thanks to Lemma 6.2 and by the definition (2.3) of $V_u(t)$, for every $t \in [0,T]$ we have $V_u(t;x) < \infty$ for $\mathcal{H}^{n-1}$-a.e. $x \in \Gamma \setminus \Gamma_N(0)$. Then, since $g(x,\xi) = \kappa(x)$ for $\xi \in [\theta(x),\infty]$, since $g(x,\cdot)$ is monotone, and since $V_u(t)$ is monotone in $t$,

$$\int_\Gamma \frac{g\big(x,V_u(t+h)\big) - g\big(x,V_u(t)\big)}{h} \, \mathrm{d}\mathcal{H}^{n-1} = \int_{\Gamma \setminus \Gamma_N(0)} \frac{g\big(x,V_u(t+h)\big) - g\big(x,V_u(t)\big)}{h} \, \mathrm{d}\mathcal{H}^{n-1}.$$

Since $V_u(t+h;x) - V_u(t;x) = \mathrm{ess}\,\mathrm{Var}([u];t,t+h)(x)$ for $\mathcal{H}^{n-1}$-a.e. $x \in \Gamma \setminus \Gamma_N(0)$ and $g'(x,V_u(t;x)) = 0$ for $\mathcal{H}^{n-1}$-a.e. $x \in \Gamma_N(0)$, by taking the limit as $h \to 0^+$ in the previous equality, by Lemma 6.2, and since $g$ is of class $\mathcal{C}^1$, we eventually deduce (6.19).

The equalities (6.18) and (6.19) combined with (EB) imply that

$$\int_{\Omega \setminus \Gamma} \nabla u(t) \cdot \nabla (\dot{u}(t) - \dot{w}(t)) \, \mathrm{d}x + \int_\Gamma g'(x,V_u(t))\big|[\dot{u}(t)]\big| \, \mathrm{d}\mathcal{H}^{n-1} = 0 \,.$$

Since $\dot{u}(t) - \dot{w}(t) = 0$ on $\partial_D \Omega$, by definition of $\partial_\nu u(t)$ we obtain

$$\int_\Gamma \partial_\nu u(t)[\dot{u}(t)] \, \mathrm{d}\mathcal{H}^{n-1} = \int_\Gamma g'(x,V_u(t))\big|[\dot{u}(t)]\big| \, \mathrm{d}\mathcal{H}^{n-1},$$

and thus

$$\int_{\{[\dot{u}(t)] \neq 0\}} \big(g'(x,V_u(t))\,\mathrm{Sign}([\dot{u}(t)]) - \partial_\nu u(t)\big)[\dot{u}(t)] \, \mathrm{d}\mathcal{H}^{n-1} = 0 \,.$$

By (iii) in Proposition 2.7, this proves the claim. $\qquad \square$

**Acknowledgements.** This material is based on work supported by the European Research Council under Grant No. 290888 "Quasistatic and Dynamic Evolution Problems in Plasticity and Fracture" and by the INdAM-GNAMPA Project 2016 "Multiscale analysis of complex systems with variational methods".

SISSA, Via Bonomea 265, 34136 Trieste, Italy
*E-mail address*, Vito Crismale: `vito.crismale@sissa.it`

*E-mail address*, Giuliano Lazzaroni: `giuliano.lazzaroni@sissa.it`

*E-mail address*, Gianluca Orlando: `gianluca.orlando@sissa.it`